\documentclass[amsfonts,onesided,11pt]{amsart}

\usepackage{amssymb,amsmath,amsfonts, amscd}
\usepackage[all,arc,cmtip]{xy}
\usepackage[bookmarks]{hyperref}
\usepackage{tikz}
\usepackage{enumerate}
\usepackage{mathrsfs}
\usepackage[margin=1in]{geometry}
\usepackage{mathabx}

\newcommand{\myind}[2]{\ensuremath{\text{ind}_{#1}^{#2}\,}}

\newcommand{\mych}{\ensuremath{\text{ch}\,}}

\newtheorem{thm}{Theorem}[subsection] 
\newtheorem{cor}[thm]{Corollary}
\newtheorem{prop}[thm]{Proposition}
\newtheorem{lem}[thm]{Lemma}
\newtheorem{conj}[thm]{Conjecture}

\theoremstyle{definition}
\newtheorem{defn}[thm]{Definition}

\newtheorem{exmp}[thm]{Example}

\newtheorem{rem}[thm]{Remark}

\makeatletter
\let\c@equation\c@thm
\makeatother
\numberwithin{equation}{subsection}

\begin{document}
\title{On Support Varieties and the Humphreys Conjecture  in type $A$}
\author{William D. Hardesty}
\address
{Department of Mathematics\\ University of Georgia \\
Athens\\ GA~30602, USA}
\email{hardesty@math.uga.edu}
\date{\today}
\begin{abstract}
Let $G$ be a reductive algebraic group scheme defined over $\mathbb{F}_p$ and let 
$G_1$ denote the Frobenius kernel of $G$. 
To each finite-dimensional $G$-module $M$, one can define the support variety $V_{G_1}(M)$, which can be
regarded as a $G$-stable 
closed subvariety 
of the nilpotent cone. 
A $G$-module is called a tilting module if it has both good and Weyl filtrations. 
 In 1997, it was conjectured by J.E. Humphreys that when $p\geq h$, the support varieties of the indecomposable tilting modules 
align with the nilpotent orbits given by the Lusztig bijection.
In this paper, we shall verify this conjecture when  $G=SL_n$ and $p > n+1$. 
\end{abstract}
\maketitle

\section{Introduction}
\subsection{}
Let $G$ be a reductive algebraic group scheme 
defined over $\mathbb{F}_p$ with Borel subgroup $B$ and maximal torus $T$. Let 
$\Phi$, $W$, $W_p=W \ltimes p\mathbb{Z}\Phi$ and 
$\mathbb{E} = \mathbb{Z}\Phi \otimes_{\mathbb{Z}}\mathbb{R}$ 
denote the root system, Weyl group, affine Weyl group and Euclidean space respectively. 
Moreover, let $G_1$ denote the Frobenius kernel and let $k$ be any algebraically closed field of characteristic $p$. 

To any finite-dimensional $G_1$-module $M$, one can associate a useful cohomological invariant called the 
\emph{support variety}, which is denoted by $V_{G_1}(M)$
(cf. \cite[Section 2.2]{npv2002} for an overview of the theory). 
 It turns out that support varieties can be identified with subvarieties of the 
\emph{$p$-restricted nullcone} 
\[
\mathcal{N}_1(G) = \{ x \in \text{Lie}(G) \mid x^{[p]} = 0\}.
\]
When $p \geq h$ (the Coxeter number of $G$), one has 
$\mathcal{N}_1(G) = \mathcal{N}(G)$, where $\mathcal{N}(G)$ is the nilpotent cone (\cite{fp}). 
So the theory of support varieties establishes a bridge between the cohomology of $G_1$-modules 
and the geometry of $\mathcal{N}_1(G)$. 

When $M$ has the structure of a $G$-module, the support variety $V_{G_1}(M)$ is 
$G$-stable. It is known that there are only finitely many $G$-orbits in $\mathcal{N}=\mathcal{N}(G)$ 
(\cite{cm}).
Hence, there are only finitely many closed subvarieties of $\mathcal{N}$ which can be realized 
as the support variety of a $G$-module.  
A major problem in representation theory has been to determine the support varieties of
various types of modules for $G$. 
Over the years, a number of results have been obtained in this direction
(cf. \cite{npv2002}, \cite{dnp}, \cite{hard}, \cite{coop2010}).
This paper will be dedicated to computing the support varieties for an important class of $G$-modules, 
known as the \emph{tilting modules} (cf.  \cite[Appendix E]{jan2003}
 for a definition and overview).
 
\subsection{}
Let $W_p^+$ denote the collection of all minimum length right coset representatives in $W \backslash W_p$. 
By introducing a certain preorder on $W_p$, one may
partition $W_p$ into \emph{two sided cells} (cf. \cite[7.15]{hum1990}).  By intersection, one also obtains a partition of $W_p^+$ into 
\emph{right cells} (cf. \cite[Theorem 1.2]{lx}).
Furthermore, there exists the \emph{Lusztig bijection}, which establishes a correspondence between the right cells of $W_p^+$ and nilpotent orbits (\cite{lusztig}).  

For $p \geq h$, it was conjectured by J.E. Humphreys in 1997 that this bijection can be realized 
by taking the
support varieties of the indecomposable tilting modules (cf. \cite[Hypothesis 12]{hum1997}). 
More precisely, 
for each $w \in W_p^+$, let $[w] \subset W_p^+$ denote the unique right cell containing $w$
and let $\mathcal{O}_{[w]} \subset \mathcal{N}$ denote the
orbit given by the Lusztig bijection. Also, for each $\lambda \in X(T)_+$, let 
$T(\lambda)$ denote the unique indecomposable tilting module with highest weight $\lambda$ 
(cf. \cite[Lemma E.3]{jan2003}). 
\begin{conj}\label{conj:hum1}
Suppose $p \geq h$, then for each $w \in W_p^+$,
$
V_{G_1}(T(w\cdot 0)) = \overline{\mathcal{O}_{[w]}}.
$
\end{conj}
\subsection{Weight cells}\label{sec:weight-cells}
For simplicity, 
assume now that $G$ is semisimple and simply connected. We will make use of the terminology to be introduced in 
Section~\ref{sec:alcove}.
It is well known that there exists a bijection between  the elements of $W_p$ (resp. $W_p^+$) and the set of  
\emph{alcoves} (resp. \emph{dominant alcoves}) of $\mathbb{E}$.
The space $\mathbb{E}$ is covered by subsets of the form $\widecheck{C}$ and 
$\mathcal{C}$ is covered by the intersections $\widecheck{C}\cap \mathcal{C}$, where
$\widecheck{C}$ denotes the lower closure of an alcove $C$. 
Thus, the right cells in $W_p^+$ can be identified
with regions in $\mathcal{C}$ called \emph{weight cells}. For each $w \in W_p^+$, they are given by 
\[
c_{[w]} = \{\lambda \in \mathcal{C} \mid \lambda \in \widecheck{y\cdot C_0},\text{ $y\in [w]$}\}.
\]

In Section~\ref{sec:translation} it will be shown that if 
$\lambda, \mu \in \widecheck{C}\cap X(T)_+$ for some alcove $C$, then 
$V_{G_1}(T(\lambda)) = V_{G_1}(T(\mu))$. 
Therefore, the following conjecture is equivalent to Conjecture~\ref{conj:hum1}. 
\begin{conj}[Humphreys Conjecture]\label{conj:hum2}
Suppose $p\geq h$ 
and $\lambda \in  c_{[w]}\cap X(T)_+$ for some $w \in W_p^+$, then
\[
V_{G_1}(T(\lambda)) = \overline{\mathcal{O}_{[w]}}.
\]
\end{conj}
In 1998,  Ostrik (\cite[Theorem 6.8]{ostrik98}) proved an analogous conjecture  for quantum groups of type $A_n$,   and in  
2006, Bezrukavnikov (\cite[3.2. Corollary 3]{bez2006}) was able to extend this result to quantum groups of any type. 
Their results are summarized in Theorem~\ref{thm:hum-quantum}. However, Conjecture~\ref{conj:hum2} still
remains open for all types, and still makes sense when $p < h$. 

%
%
\subsection{}
One of the major obstacles to proving Conjecture~\ref{conj:hum2}
 has been the difficulty of determining the weight cells
$c_{[w]}$.
Although when $G$ is a reductive group of type $A_n$ (i.e., $G=SL_{n+1}(k)$), there is a result to due to 
Shi which gives 
an explicit description of the weight cells; it is described in Section~\ref{sec:cell_regions} (cf. \cite{shi} for the original result).
In fact, progress has already been made in the type $A_n$ situation by Cooper, 
who first extended Conjecture~\ref{conj:hum2} by removing the assumption that $p\geq h$.
Cooper then made significant progress in verifying this conjecture for small primes, including a complete verification when $p=2$
(cf. \cite[Theorem 7.3.1]{coop2010}). 

To be more specific, in the type $A_n$  case, it is well-known that the orbits of $\mathcal{N}$ correspond to partitions $\pi \vdash n+1$  (\cite{cm}).
So the Lusztig bijection establishes a correspondence between weight cells and partitions. 
In Definition~\ref{def:hum-part} and Remark~\ref{rem:weight-cell}, these weight cells will be explicitly 
described by associating a partition, $s(\lambda)$, to each $\lambda \in \mathcal{C}$.
The main result of this paper is the following theorem, which verifies Conjecture~\ref{conj:hum2} in type $A_n$
for all $n$, when $p>h=n+1$. 
\begin{thm}\label{thm:main}
 Let $G=SL_{n+1}(k)$ with $p > n+1$, then  for each $\lambda \in X(T)_+$, 
 \[
 V_{G_1}(T(\lambda)) =\overline{\mathcal{O}_{s(\lambda)^t}}.
 \]
\end{thm}

The proof of Theorem~\ref{thm:main} will  begin by showing that 
$V_{G_1}(T(\lambda)) \subseteq \overline{\mathcal{O}_{s(\lambda)^t}}$ for each $\lambda \in X(T)_+$, 
which places an upper bound to the support variety $V_{G_1}(T(\lambda))$. 
This will require several steps: first in Section~\ref{sec:alcove} some
results regarding the alcove geometry associated to the affine Weyl group will be obtained and 
we will define the weak order on alcoves (see Definition~\ref{defn:super-linkage}). 
By recalling a few key
identities involving translation functors and wall crossing functors
in Section~\ref{sec:translation}, Corollary~\ref{cor:sup-cor} will relate
this order relation to the ordering of support varieties by inclusion. 
All of the results in these two sections will hold for arbitrary simple, simply connected groups. 

In Section~\ref{sec:cell_regions}, an explicit description
of the weight cells in type $A_n$ will be presented. 
The main result of this section is
Proposition~\ref{cor:good-part}, which gives
an equivalent characterization of the partitions $s(\lambda)$. 
Finally, Section~\ref{sec:upper-bound} will include a proof of Proposition~\ref{prop:upperbound}, which
establishes the upper bound portion of Theorem~\ref{thm:main} under the slightly relaxed assumption
that $p\geq n+1$. 

The remainder of the paper will be dedicated to establishing the lower bound. 
Section~\ref{sec:quantum} will review the necessary definitions and facts about quantum groups. 
It will include a result by Andersen, which allows one to ``lift'' tilting modules over $G$ to tilting modules over an analogous 
 quantum group  (cf. \cite[5.3]{and1998}).
Some results and conjectures regarding the 
support varieties of these lifted tilting modules will also be presented, 
including a complete description of them in type $A$ 
(see Proposition~\ref{prop:quantum-upper}).
Section~\ref{sec:type-A} will make use of the fact that in type $A$, every non-dense nilpotent orbit intersects a proper Levi factor
(see Lemma~\ref{lem:a_n-levi}). Proposition~\ref{prop:lowerbound} will give the lower bound, which will follow by making direct comparisons to the quantum case when $p>n+1$. This proposition, along with Proposition~\ref{prop:upperbound}, will yield Theorem~\ref{thm:main}.
\begin{rem}
It is useful to note that the support variety calculations for induced modules in \cite{npv2002}, the irreducible modules in \cite{dnp}
and the higher sheaf cohomology modules in \cite{hard} made explicit use of the known character formulas for the corresponding modules, to get 
the lower bound.  However, character formulas for the indecomposable tilting modules, when $G$ is not of type $A_1$, have yet to be determined.
In fact, to the author's best knowledge, there is no known conjecture which predicts the characters of all indecomposable tilting modules for arbitrary semisimple groups (see \cite{lw}  and \cite[3.6: Remark (i)]{and1997} for some partial results and conjectures). The lower bound calculation given in Proposition~\ref{prop:lowerbound}, will only utilize partial information about the characters. Namely, it will use the character
formulas given by Soergel in \cite{soergel1} and \cite{soergel2} for quantum groups, and the identity \eqref{eqn:lifting}. 
\end{rem}
\subsection{Acknowledgements} 
This paper is a part of the author's PhD dissertation and he would like thank
his PhD thesis advisor, Daniel Nakano, for all of his consultation during this project. I would also like to thank Henning Haahr Andersen, 
Jim Humphreys and William Graham, for their very helpful feedback. 
The author was partially supported by the Research Training Grant, DMS-1344994, from the NSF.

\section{Alcove geometry}\label{sec:alcove}
\subsection{}
In this section, 
 assume that $\Phi$ is an irreducible root system of any type and that $p\geq 1$ is any positive integer. 
Let $\Delta = \{\alpha_1,\dots, \alpha_n\}$ and 
$
\Phi^+ = \mathbb{N}\Delta \cap \Phi
$
 denote the basis and 
the set of positive roots respectively. 
Let $\mathbb{E} = \mathbb{Z}\Phi \otimes_{\mathbb{Z}} \mathbb{R}$ be the Euclidean space.
Then $\mathbb{E}$ is given the lattice ordering, where
for
$\lambda, \mu \in \mathbb{E}$,
$\lambda \leq \mu$ will be taken to mean that $\mu - \lambda \in \mathbb{N}\Phi^+$.
 Let $\alpha_0 \in \Phi^+$ denote the maximal short 
root with respect to this ordering. 
The \emph{strong linkage relation} also gives an ordering on $\mathbb{E}$, where
$\lambda \uparrow \mu$ will denote when $\lambda$ is \emph{strongly linked} to $\mu$ (cf. \cite[II.6]{jan2003}). 


The affine Weyl group $W_p = W \ltimes p\mathbb{Z}\Phi$ is a Coxeter group with generators 
$\mathcal{S}=\{s_0,s_1,\dots,s_n\},$ where $s_0$ denotes the \emph{affine reflection} and the generators $s_1,\dots s_n$ correspond
to the basis elements $\alpha_1,\dots,\alpha_n$ (see Definition~\ref{defn:walls} below). The group $W_p$ is equipped with the standard 
\emph{length function} $\ell: W_p \rightarrow \mathbb{N}$, where $\ell(w)$ denotes the length of any reduced expression for $w \in W_p$. 
Moreover, $W_p$ acts on $\mathbb{E}$ 
 by both the linear action 
 and the \emph{dot action}. 
As usual, the linear action will be denoted by $\lambda \mapsto w(\lambda)$ and the dot action will be denoted by
$\lambda \mapsto w\cdot \lambda = w(\lambda + \rho)-\rho$, where $\rho \in \mathbb{E}$ is the half sum of
the positive roots. 

 The group $W_p$ is partially ordered by
the Bruhat ordering, which will be denoted by $\leq$. 
 Let $W_p^+$ be the 
set of minimal length right cosets for the finite Weyl group $W$ in $W_p$, and let 
\[
\mathcal{C} = \{ \lambda \in \mathbb{E} \mid \langle \lambda + \rho, \alpha^{\vee}\rangle >0 \text{ for all $\alpha \in \Delta$}\}
\] 
denote the \emph{dominant chamber} of 
$\mathbb{E}$.
For each $\alpha \in \Phi^+$ and $n \in \mathbb{Z}$, let
\[
  H_{\alpha,np} = \{ \lambda \in \mathbb{E} \mid \langle \lambda + \rho, \alpha^{\vee} \rangle = np \},
\]
let $s_{\alpha,np} \in W_p$ denote the affine reflection across $H_{\alpha,np}$ (with respect to the dot action),
and let $ \mathcal{H}=\bigcup_{\alpha \in \Phi^+,\, n \in \mathbb{Z}}H_{\alpha,np}$.
  The connected components of $\mathbb{E}\backslash \mathcal{H}$ 
are called \emph{alcoves}. Let $\mathcal{A}$ denote the collection of all the alcoves. The collection of 
\emph{dominant alcoves} will be denoted by $\mathcal{A}^+$; it consists of all the alcoves contained in 
 $\mathcal{C}$. 
 The dot action by $W_p$ on $\mathbb{E}$ induces a simply transitive action on $\mathcal{A}$ by 
 sending $C \mapsto w\cdot C$ for any $C \in \mathcal{A}$ and $w \in W_p$. 
The set  $\mathcal{A}$ also has an ordering induced by  the strong linkage relation on $\mathbb{E}$, where
 $C_1 \uparrow C_2$ if there exists $\lambda_i \in  C_i$ such that $\lambda_1 \uparrow \lambda_2$ for $i=1,2$. 
 
 If $C$ is any alcove, then it is uniquely defined by a set of  integers $\{n_{\alpha}\}_{\alpha \in \Phi^+}$, where
\[
  C = \{ \lambda \in \mathbb{E} \mid  (n_{\alpha}-1)p < \langle \lambda + \rho, \alpha^{\vee} \rangle < n_{\alpha}p,
  \, n_{\alpha} \in \mathbb{Z}, \, \alpha \in \Phi^+ \}.
\]
The \emph{upper closure} of $C$ is given by
\[
  \widehat{C} = \{ \lambda \in \mathbb{E} \mid  (n_{\alpha}-1)p < \langle \lambda + \rho, \alpha^{\vee} \rangle \leq n_{\alpha}p,
  \, n_{\alpha} \in \mathbb{Z}, \, \alpha \in \Phi^+ \}.
\]
The \emph{lower closure}  of $C$ is defined to be
\[
  \widecheck{C} = \{ \lambda \in \mathbb{E} \mid  (n_{\alpha}-1)p \leq \langle \lambda + \rho, \alpha^{\vee} \rangle < n_{\alpha}p,
  \, n_{\alpha} \in \mathbb{Z}, \, \alpha \in \Phi^+ \}.
\]
For each $\lambda \in \mathbb{E}$, the unique alcove satisfying 
$\lambda \in \widecheck{C(\lambda)}$ is denoted by $C(\lambda)$. 
The alcove given by $n_{\alpha}=1$ for all $\alpha \in \Phi^+$ is called the \emph{bottom alcove} and will
be denoted by $C_0$. 
\begin{rem}\label{rem:bijection}
One obtains a bijection between $\mathcal{A}$ (resp. $\mathcal{A}^+$) and 
 $W_p$ (resp. $W_p^+$) by identifying $w \leftrightarrow w\cdot C_0$. 
  This bijection identifies the  strong linkage relation on $\mathcal{A}^+$ with the Bruhat order on $W_p^+$, where $w_1\cdot C_0 \uparrow w_2\cdot C_0$ if and only if 
 $w_1\leq w_2$ (cf. \cite[C.1]{jan2003}). 

 \end{rem}
 

If $\{n_{\alpha}\}_{\alpha \in \Phi^+}$ is a set of integers defining an alcove as above, then
for any subset $S\subseteq \Phi^+$,
\[
 F =  \left\{ \lambda \in \overline{C}\,\, \middle| \,\,
 	  \begin{aligned}
	    \langle \lambda + \rho, \alpha^{\vee} \rangle &= n_{\alpha}p, \,\,\text{ if  $\alpha \in S$} \\
	    (n_{\alpha}-1)p < \langle \lambda + \rho, \alpha^{\vee} \rangle &< n_{\alpha}p, \,\, \text{ if  $\alpha \in \Phi^+\backslash S$}
	     \end{aligned} \right\}
\]
is called a \emph{facette}, where $F\subset \widehat{C}$. 
We can similarly define the lower closure (resp. upper closure) $\widecheck{F} \subseteq \overline{F}$ 
(resp. $\widehat{F}\subseteq \overline{F}$).   The collection of all facettes in $\mathbb{E}$ will be denoted by $\mathcal{F}$, 
where $\mathcal{F}$ is also acted on by $W_p$ via $F \mapsto w\cdot F$ for $w \in W_p$ and $F \in \mathcal{F}$.  
Each alcove is a facette in its own right, and thus
$\mathcal{A} \subseteq \mathcal{F}$. 

For any $\lambda \in \mathbb{E}$, denote the unique facette containing $\lambda$ by 
$F(\lambda)$. Every non-empty facette 
 is of the form $F=F(\lambda)$ for some $\lambda \in \mathbb{E}$. 
Also, let 
\[
\operatorname{Stab}_{W_p}(\lambda) = \{ w \in W_p \mid w\cdot \lambda = \lambda\}
\]
 denote the 
\emph{stabilizer subgroup} of $\lambda \in \mathbb{E}$. If the facette $F=F(\lambda)$ 
is given by $S\subseteq \Phi^+$ and $\{n_{\alpha}\}_{\alpha \in \Phi^+}$ as above, then 
$\operatorname{Stab}_{W_p}(\lambda)$ is generated by the set of reflections $\{s_{\alpha,n_{\alpha}p}\}_{\alpha \in S}$.
 \begin{defn}\label{defn:walls}
The \emph{walls} of $C_0$ are defined to be the hyperplanes: $H_0 = H_{\alpha_0,p}$ (called the 
\emph{affine wall}) and $H_i = H_{\alpha, 0}$ for $\alpha \in \Phi^+$. 
 The elements of $\mathcal{S}$ are the reflections across the walls of $C_0$, which are given by $s_0 = s_{\alpha_0, p}$ and $s_i = s_{\alpha_i,0}$ for $i=1,\dots,n$.
More generally, for an alcove $C = w\cdot C_0$ with $w \in W_p$, the \emph{walls} of $C$ are the hyperplanes 
$w\cdot H_{i}$ for $i=0,\dots, n$. 
\end{defn}

\begin{rem}\label{rem:walls}
If $C$ is an alcove defined by the integers $\{n_{\alpha}\}_{\alpha \in \Phi^+}$, then a hyperplane $H_{\alpha,mp}$
is a wall of $C$ (see Definition~\ref{defn:walls}), if there is an element $\lambda \in H_{\alpha,mp}\cap \overline{C}$ satisfying 
$\operatorname{Stab}_{W_p}(\lambda) = \{1, s_{\alpha,mp}\}$. 
It follows that $m \in \{n_{\alpha}-1,n_{\alpha}\}$, 
where if $m= n_{\alpha}$ (resp. $m=n_{\alpha}-1$), then $H_{\alpha, mp}$ is called an \emph{upper wall} 
(resp. \emph{lower wall}) of $C$.  Moreover, $C\uparrow s_{\alpha,mp}\cdot C$ if and only if 
$H_{\alpha,mp}$ is an upper wall of $C$.  
\end{rem}
%

\begin{lem}\label{lem:lclosure}
For each facette $F=F(\lambda)$, 
\[
\widecheck{F} = \{ \mu \in \overline{F} \mid \lambda \geq w\cdot \lambda\, \text{ for all 
									$w \in \operatorname{Stab}_{W_p}(\mu)$}\}.
\]
\end{lem}
\begin{proof}
Let $\widetilde{F}$ denote the right hand side of the stated identity. Suppose
$F$ is given by the data $S\subseteq \Phi^+$ and $\{n_{\alpha}\}_{\alpha \in \Phi^+}$. 
In other words, for all $\alpha \not\in S$,
$(n_{\alpha}-1)p < \langle \lambda + \rho, \alpha^{\vee} \rangle < n_{\alpha}p$
and for all $\alpha \in S$, $\langle \lambda + \rho, \alpha^{\vee} \rangle = n_{\alpha}p$.

 Let us first 
show $\widetilde{F} \subseteq \widecheck{F}$. By definition 
$\widetilde{F} \subseteq \overline{F}$,  so it suffices to show
that $\widetilde{F}$ doesn't contain any elements of the form $\mu \in \overline{F}$ satisfying
$\langle \mu + \rho, \alpha^{\vee} \rangle = n_{\beta}p$
for some $\beta \not\in S$. However, if such a $\mu$ exists, then the reflection 
$s_{\beta,n_{\beta}p} \in \text{Stab}_{W_p}(\mu)$ and $\lambda \nleq s_{\beta,n_{\beta}p}\cdot \lambda$,
since
$s_{\beta,n_{\beta}p}\cdot \lambda = \lambda + m\beta$,
where $m = n_{\beta}p-\langle \lambda + \rho, \beta^{\vee} \rangle > 0$.

To prove $\widetilde{F}\supseteq \widecheck{F}$, begin by choosing an arbitrary element $\mu \in \widecheck{F}$.
Let $T\subseteq \Phi^+$ be the collection of all $\beta \in \Phi^+$ satisfying 
$\langle \mu + \rho, \beta^{\vee} \rangle = (n_{\beta}-1)p$.
Then the stabilizer of $\mu$ is the subgroup of $W_p$ generated by the reflections 
$s_{\alpha, n_{\alpha}p}$ for $\alpha \in S$ and $s_{\beta,(n_{\beta}-1)p}$ for $\beta \in T$. 
However, for all $\alpha \in S$, $s_{\alpha,n_{\alpha}p}\cdot \lambda = \lambda$ and for all $\beta \in T$,
$s_{\beta,(n_{\beta}-1)p}\cdot \lambda = \lambda + m\beta \leq \lambda$,
since $m = (n_{\beta}-1)p-\langle \lambda + \rho, \beta^{\vee} \rangle < 0$. Therefore,
$\mu \in \widetilde{F}$.
\end{proof}
If $F$ is a  facette, there exists a unique alcove $C$ such that $F\subseteq \widecheck{C}$. Moreover, 
if $F$ is given by the data $S\subseteq \Phi^+$ and $\{n_{\alpha}\}_{\alpha \in \Phi^+}$, 
then $C$ is given by the integers $\{m_{\alpha}\}_{\alpha \in \Phi^+}$, where $m_{\alpha} = n_{\alpha}$ if $\alpha \not\in S$ and 
$m_{\alpha} = n_{\alpha} +1$ if $\alpha \in S$.  Concretely, 
\[
 C =  \left\{ \lambda \in \mathbb{E}\,\, \middle| \,\,
 	  \begin{aligned}
	    n_{\alpha}p< \langle \lambda + \rho, \alpha^{\vee} \rangle &< (n_{\alpha}+1)p, \,\,\text{ if  $\alpha \in S$} \\
	    (n_{\alpha}-1)p < \langle \lambda + \rho, \alpha^{\vee} \rangle &< n_{\alpha}p, \quad\quad\,\,\,\,\, \text{ if  $\alpha \in \Phi^+\backslash S$}
	     \end{aligned} \right\}.
\]
\subsection{}
 A useful refinmenent of the Bruhat ordering,  called the \emph{weak ordering}, can be placed on $W_p$ (see \cite[5.9]{hum1990}). 
 \begin{defn}\label{defn:super-linkage}
 Let $w_1, w_2 \in W_p$ be arbitrary and suppose $w_1 = t_1\cdots t_{m_1}$ is a reduced expression, then $w_1 \preceq w_2$ if and only if there is a sequence of elements
  \[
  w_1 = w'_0 \leq w'_1 \leq \cdots \leq w'_{m_2} = w_2,
  \]
  where $w'_i \in W_p$  and $w'_i = t_1\cdots t_{m_1+i}$ is a reduced expression for $i=0,\dots, m_2$. 
%
 For any two alcoves $C_1, C_2 \in \mathcal{A}$, take $C_1 \preceq C_2$ to mean that 
 $w_1\preceq w_2$ for the unique $w_1,w_2 \in W_p$ satisfying $w_i \cdot C_0 = C_i$ for $i=1,2$. 
This defines the weak order on $W_p$ (respectively $\mathcal{A}$). 
 \end{defn}
 The weak ordering restricts to give an order relation on $W_p^+$ and $\mathcal{A}^+$.
%
%

\begin{rem}\label{rem:super-linkage}
   If $C \in \mathcal{A}$ is of the form $C=w\cdot C_0$ for some $w \in W_p$, then for each $i=0,\dots, n$ and 
   $s_i \in \mathcal{S}$,  
   the alcove $ws_{i}\cdot C$ is obtained by reflecting $C$ across the wall $w\cdot H_i$ (see Definition~\ref{defn:walls}).
   Furthermore, if $w, ws_i \in W_p^+$, then $w \leq ws_i$ if and only if $w\cdot H_i$ is an upper wall of $C$
   (see Remark~\ref{rem:walls}). 
   Hence, if $C_1=w_1\cdot C_0$ and $C_2 = w_2\cdot C_0$ are two dominant alcoves, then 
   $C_1 \preceq C_2$ if and only if there is a sequence of alcoves 
   \[
   C_1 = C'_0\uparrow  C'_1\uparrow \dots \uparrow C'_{m_2} = C_2,
   \]
   where $C'_i = t_1\cdots t_{m_1+i}\cdot C_0$ for $i=0,\dots, m_2$. 
%
   Thus, by Remark~\ref{rem:walls}, for each $i=0,\dots, m_2-1$, 
    $C'_{i+1} = s_{\beta_i, n_ip}\cdot C'_i$, where $\beta_i \in \Phi^+$ and 
   $H_{\beta_i, n_ip}$ is an upper wall of $C'_i$.  
   \end{rem}

The following lemma gives an important 
characterization of the weak order.
\begin{lem}\label{lem:ss-implies-s}
  Let $C_1,C_2 \in \mathcal{A}^+$  be two alcoves defined by the nonnegative integers $\{n_{\alpha}\}_{\alpha \in \Phi^+}$ and 
  $\{m_{\alpha}\}_{\alpha \in \Phi^+}$ respectively.
 Then $C_1\preceq C_2$ if and only if $n_{\alpha} \leq m_{\alpha}$ for all $\alpha \in \Phi^+$. 
\end{lem}
\begin{proof}
Let $C_1 = w_1 \cdot C_0$ and $C_2 = w_2 \cdot C_0$, and suppose $C_1\preceq C_2$. Thus,
$w_1 = t_1\cdots t_{m_1}$ and 
$w_2 = t_1 \cdots t_{m_1+m_2}$ as in Definition~\ref{defn:super-linkage}. 
Let $C'_i = t_1\cdots t_{m_1+i}\cdot C_0 \in \mathcal{A}^+$ for $i=0,\cdots, m_2$.
For each $i$, the alcove $C'_i$  is defined by the set of integers
$\{(n_{\alpha})_i\}_{\alpha\in\Phi^+}$. 
By Remark~\ref{rem:super-linkage},  there exists a root $\beta_i \in \Phi^+$ such that $C'_{i+1}=s_{\beta_i, (n_{\beta_i})_ip}\cdot C'_i$, where $H_{\beta_i, (n_{\beta_i})_ip}$ is an upper wall of $C'_{i}$. Thus,
\[
(n_{\alpha})_{i+1} = \begin{cases}
	         (n_{\alpha})_i & \text{ if  $\alpha \neq \beta_i$} \\
	         (n_{\alpha})_i+1  & \text{ if $\alpha = \beta_i$}.
	    \end{cases}
\]
 It follows that $n_{\alpha}\leq (n_{\alpha})_1 \leq \cdots \leq (n_{\alpha})_{m_2}= m_{\alpha}$ for all $\alpha \in \Phi^+$. 

For the converse,  perform induction on 
$d = \sum_{\alpha \in \Phi^+} m_{\alpha} - n_{\alpha}$. Observe that since $m_{\alpha}\geq n_{\alpha}$ for all 
$\alpha \in \Phi^+$,
then $d$ is equal to the number of hyperplanes separating $C_1$ and $C_2$. This is because the complete set of hyperplanes separating $C_1$ and $C_2$ is given by 
\begin{equation}\label{eqn:seperation}
\{H_{\alpha, kp} \mid \alpha \in \Phi^+, \, n_{\alpha} \leq k \leq m_{\alpha}-1\},
\end{equation}
and has $d$ elements. 

 The case where $d=0$ holds because $C_1 = C_2$ implies $C_1 \preceq C_2$. 
Now for the inductive step,  suppose $d \geq 1$. Then $C_1\neq C_2$ since $d > 0$, and so
there must exist some wall of $C_1$ which separates $C_1$ and $C_2$. 
If this wall is given by $H_{\beta,mp}$ for some  $\beta \in \Phi^+$, then by Remark~\ref{rem:walls},  
$m \in \{n_{\beta}-1,n_{\beta}\}$, and by \eqref{eqn:seperation},
 $n_{\beta} \leq m \leq m_{\beta}-1$. Thus, $m = n_{\beta}$, and hence $H_{\beta, mp}$ is an upper wall of $C_1$.
 The alcove given by $C = s_{\beta, mp}\cdot C_1 \in \mathcal{A}^+$ 
  is defined by the integers $\{r_{\alpha}\}_{\alpha \in \Phi^+}$, 
where
\[
r_{\alpha} = \begin{cases}
	         n_{\alpha} & \text{ if  $\alpha \neq \beta$} \\
	         n_{\alpha}+1  & \text{ if $\alpha = \beta$}.
	    \end{cases}
\]
Thus, $r_{\alpha} \leq m_{\alpha}$ for all $\alpha \in \Phi^+$, $\sum_{\alpha \in \Phi^+} m_{\alpha}-r_{\alpha} = d-1$ and,
by Remark~\ref{rem:super-linkage}, $C_1 \preceq C$. 
By the inductive hypothesis, $C \preceq C_2$, and therefore 
$C_1 \preceq C_2$.  
%
\end{proof}
By Remark~\ref{rem:bijection}, it follows that if
$C_1 \preceq C_2$, then $C_1 \uparrow C_2$. Since by definition, if $w_1\preceq w_2$,  then $w_1 \leq w_2$. 
   However, $C_1\uparrow C_2$ doesn't generally imply $C_1\preceq C_2$ 
   (see the argument in  Remark~\ref{rem:link-counter} for a counterexample). 
\subsection{Stabilizer subroot systems}
Each $\lambda \in \mathbb{E}$ can be associated to a certain subroot system of $\Phi$.  
\begin{defn}\label{defn:stab-ss}
For each $\alpha \in \Phi$, let 
$d_{\alpha} = \langle \alpha, \alpha\rangle/ \langle \alpha_0, \alpha_0\rangle \in \{1,2,3\}$. 
For $\lambda \in \mathbb{E}$, define
\[
\Phi_{\lambda,p} = \{ \alpha \in \Phi \mid d_{\alpha}\,\langle \lambda + \rho, \alpha^{\vee}\rangle \in p\mathbb{Z}
\}
\]
to be the \emph{stabilizer subroot system of $\lambda$}. 
\end{defn}
Now $\Phi_{\lambda,p}$ is a closed subroot system of $\Phi$ because
 for any two roots $\alpha, \beta \in \Phi$,  
 $d_{\alpha + \beta}(\alpha + \beta)^{\vee} = d_{\alpha}\alpha^{\vee} + d_{\beta}\beta^{\vee}$ (cf. \cite[6.2]{npv2002}).
 Subroot systems can also be associated to facettes, since if $F \in \mathcal{F}$ is a facette and $\lambda, \mu \in F$, 
 then 
$\Phi_{\lambda,p} = \Phi_{\mu,p}$. 
\begin{defn}
For each $F \in \mathcal{F}$ and any $\lambda \in F$, the  \emph{stabilizer subroot system of $F$} is given by
 $\Phi(F) = \Phi_{\lambda,p}$.
\end{defn}

Let $X \subset \mathbb{E}$ denote the lattice consisting of all $\lambda \in \mathbb{E}$ satisfying 
$\langle \lambda + \rho, \alpha^{\vee}\rangle \in \mathbb{Z}$ for any $\alpha \in \Phi$, then
$\mathbb{Z}\Phi \subset X \subset \mathbb{E}$. 
A basis of $X$ is given by the \emph{fundamental weights} $\omega_1,\dots, \omega_n \in X$,
which satisfy $\langle \omega_i, \alpha_j^{\vee}\rangle =\delta_{ij}$.  
In fact, if $\Phi$ is the root system associated to a semisimple, simply connected algebraic group $G$, then 
$X = X(T)$ is the weight lattice for $G$. 
The \emph{extended affine Weyl group} is given by $\widetilde{W_p} = W \ltimes pX$, where  $W_p \unlhd \widetilde{W_p}$.
The dot action of $\widetilde{W_p}$ on $\mathbb{E}$ induces an action on $\mathcal{F}$, where for each 
$F \in \mathcal{F}$
\begin{equation}\label{eqn:extended_orbits}
\widetilde{W_p}\cdot F = \{F' \in \mathcal{F} \mid \Phi(F') = w(\Phi(F)) \text{ for some $w \in W$}\}. 
\end{equation}
\subsection{Lattice points}\label{ssec:lattice}
In this subsection, we assume that $p$ is a prime number. 
\begin{defn}\label{defn:good-prime}
A prime $p$ is said to be \emph{good} for a root system $\Phi$ if 
for any closed subroot system $\Phi' \subseteq \Phi$, the quotient $\mathbb{Z}\Phi/ \mathbb{Z}\Phi'$ has no 
$p$-torsion. Equivalently, $p$ is good unless $\Phi$ has a component of type $B_n$, $C_n$, $D_n$ and $p=2$; 
$\Phi$ has a component of type $E_6$, $E_7$, $F_4$, $G_2$ and $p=2, 3$; or $\Phi$ has a component of type $E_8$ and $p=2,3,5$. 
\end{defn}
\begin{rem}
It follows that $p$ is good for $\Phi$ if and only if $p$ is good for $\Phi^{\vee}$. 
Also, if $p \geq h$, then $p$ is a \emph{good} prime for $\Phi$ (and equivalently for $\Phi^{\vee}$).
\end{rem}
It will be useful to determine precise conditions under which a  facette 
$F$  satisfies $F\cap X \neq \emptyset$, when $p$ is good (or $p \geq h$).
Notice that if $\lambda \in X$, then 
  $\langle \lambda + \rho, \alpha^{\vee} \rangle \in \mathbb{Z}$ for all $\alpha \in \Phi$, and thus
  $\Phi_{\lambda,p}^{\vee}$ is also a closed subroot system of $\Phi^{\vee}$(cf. \cite[6.2]{npv2002}).
  Furthermore, since $p$ is good and $\mathbb{Z}\Phi^{\vee}/\mathbb{Z}\Phi_{\lambda,p}^{\vee}$ contains no $p$-torsion, 
  it can verified that
   \begin{equation}\label{en:good-ss}
 \Phi_{\lambda,p} = \{\alpha \in \Phi \mid \langle \lambda + \rho, \alpha^{\vee}\rangle \in p\mathbb{Z} \}.
 \end{equation}
  \begin{defn}\label{defn:parabolic}
  A \emph{parabolic subroot system} of $\Phi$ is defined to be a subroot system of the form
\[
\Phi_I = \mathbb{Z}I\cap \Phi,
\]
where $I\subseteq \Delta$ is any subset. 
\end{defn}

The following lemma gives necessary conditions for when $F\cap X \neq \emptyset$. 
\begin{lem}\label{lem:nec-lattice}
Let $p$ be good and suppose $F \in \mathcal{F}$ satisfies $F\cap X \neq \emptyset$, then 
$\Phi(F) = w(\Phi_I)$ for some $w \in W$ and some $I \subseteq \Delta$. 
\end{lem}
\begin{proof}
Let $\lambda \in F\cap X$, then $\Phi(F) = \Phi_{\lambda,p}$. By the comments
above,  $\Phi^{\vee}_{\lambda,p}$ is a closed subroot 
system of $\Phi^{\vee}$ and $\mathbb{Z}\Phi^{\vee}/\mathbb{Z}\Phi_{\lambda,p}^{\vee}$ contains no $p$-torsion. 
If $\alpha \in \Phi$ 
 satisfies $m\alpha^{\vee} \in \Phi_{\lambda, \alpha}^{\vee}$ for some $m \in \mathbb{Z}$, then 
 \[
 m\,\langle \lambda + \rho, \alpha^{\vee} \rangle = \langle \lambda + \rho, m\alpha^{\vee} \rangle \in p\mathbb{Z}.
 \]
 Since $\mathbb{Z}\Phi^{\vee}/\mathbb{Z}\Phi_{\lambda,p}^{\vee}$ contains no $p$-torsion, then  $p \nmid m$.
 It follows that  
  $\alpha^{\vee} \in \Phi_{\lambda,p}^{\vee}$, and thus   $\alpha \in \Phi_{\lambda,p}$. Therefore, $\Phi_{\lambda,p} = \mathbb{Q}\Phi_{\lambda,p}\cap \Phi$ and, by \cite[Proposition 24, p. 165]{bourbaki}, 
 there exists $I \subseteq \Delta$ and $w \in W$ such that $\Phi_{\lambda,p} = w(\Phi_I)$. 
\end{proof}

When $p\geq h$, necessary and sufficient conditions for when $F\cap X \neq \emptyset$ 
can be determined. 
\begin{prop}\label{prop:latticept}
Let $p \geq h$ be prime, then $F \in \mathcal{F}$ satisfies $F\cap X \neq \emptyset$ if and only if 
$\Phi(F) = w(\Phi_I)$ for some $w \in W$. 
\end{prop}
 \begin{proof}
 First suppose that $\Phi(F) = w(\Phi_I)$ for some $w \in W$ and $I\subseteq \Delta$. If
 $\lambda + \rho = \sum_{\alpha_i \not\in I} \omega_i$, then since $p\geq h$,
 $\Phi_{\lambda,p} = \Phi_I$. By \eqref{eqn:extended_orbits}, there exists $x \in \widetilde{W_p}$ such that 
 $F =x\cdot F(\lambda)$. Then $x\cdot \lambda \in F\cap X$. 
 
 Since $p$ is good, the converse follows from Lemma~\ref{lem:nec-lattice}.
 
 \end{proof}

 \section{Translation functors and tensor ideals }\label{sec:translation}
\subsection{}
In this section, translation functors will be employed to establish some identities relating $\preceq$ to 
 inclusions of thick tensor ideals and support varieties of tilting modules (see the definition below). 

\begin{defn}\label{def:tilting-category}
For any reductive group $G$ over a field $k$ of characteristic $p>0$, let
 $\mathcal{T}=\mathcal{T}(G)$ denote the full subcategory in the category of rational $G$-modules, consisting of all 
 finite-dimensional tilting modules for $G$.

For each $M \in \mathcal{T}$, the \emph{thick tensor ideal} generated by $M$ is given by
\[
 \langle M \rangle = \{ N \in \mathcal{T} \mid  N \, | \,M\otimes L\, \text{ for some $L \in \mathcal{T}$}\}.
 \]
 Where $M \mid N$ for $M, N \in \mathcal{T}$, denotes the 
 existence of a decomposition of the form $N = M \oplus T$ for some $T \in \mathcal{T}$.  
 \end{defn}
  
  If $M_1, M_2 \in \mathcal{T}$
  satisfy
 $\langle M_1 \rangle \subseteq \langle M_2 \rangle$, then $V_{G_1}(M_1) \subseteq V_{G_1}(M_2)$.  Now recall
 the definition of 
 translation functors (see also Lemma~\ref{lem:sup-ideal}).
\begin{defn}
Let $M$ be a rational $G$-module. Then for any $\lambda, \mu \in X(T)_+ \cap \overline{C}$ for some alcove $C$, 
\[
 T_{\lambda}^{\mu}M = \text{pr}_{\mu}(M \otimes T(\nu)),
 \]
 where $\text{pr}_{\mu}$ denotes the projection onto the block containing the simple module $L(\mu)$.
  If $\lambda_0, \mu_0 \in \overline{C_0}$ with 
  $\lambda_0 \in W_p\cdot \lambda$ and $\mu_0\in W_p \cdot \mu$, then 
 $\nu \in X(T)_+$ satisfies $\nu = w(\mu_0-\lambda_0)$ for some $w \in W$ 
 (cf. \cite[II.7: Definition 7.6 and the remarks following Lemma 7.7]{jan2003}).
\end{defn}

The following proposition is a clarification of \cite[Proposition E.11]{jan2003}, which contains a minor mistake in the statement.  
It provides some information about the behavior of tilting modules under translation functors. 
 \begin{prop}\label{prop:translation1}
Let $\lambda, \mu \in X(T)_+$ satisfy $\mu \in \widecheck{F(\lambda)}$, 
 then $T_{\mu}^{\lambda}T(\mu) \cong T(\lambda)$
  and 
  \begin{equation}\label{eqn:trans}
  T_{\lambda}^{\mu}T(\lambda) = [\operatorname{Stab}_{W_p}(\mu):\operatorname{Stab}_{W_p}(\lambda)]\,T(\mu).
  \end{equation}
  In particular, 
  $\langle T(\lambda)\rangle= \langle T(\mu)\rangle$.
\end{prop}
\begin{proof}
By Lemma~\ref{lem:lclosure}, $\lambda \geq w\cdot \lambda$ for all $w \in \operatorname{Stab}_{W_p}(\mu)$. The proof 
of  \cite[Proposition E.11]{jan2003} then implies \eqref{eqn:trans}. 
\end{proof}

In general, if $\mu \in \overline{F(\lambda)}\backslash \widecheck{F(\lambda)}$, then the structure of modules
such as $T_{\lambda}^{\mu}T(\lambda)$ or 
$T_{\mu}^{\lambda}T_{\lambda}^{\mu}T(\lambda)$ is much more difficult to describe. However, the following lemma
gives some insight into this case. 

\begin{prop}\label{lem:wall_crossing}
  Let $\lambda, \mu \in X(T)_+$ and suppose $\mu \in \overline{F(\lambda)}$.
   If $\lambda'$ is the maximal element of $\text{Stab}_{W_p}(\mu)\cdot \lambda$,
  then it is the highest weight of the tilting module
  $\Theta\, T(\lambda) = T_{\mu}^{\lambda}T_{\lambda}^{\mu}T(\lambda)$, and thus 
 $T(\lambda') \mid \Theta\, T(\lambda)$. 
\end{prop}
\begin{proof}
  By definition, $T(\lambda )$ has a filtration 
  \[
    0=F_0 \subset \cdots \subset F_m=T(\lambda)
  \]
  where $F_i/F_{i-1} \cong V(x_i\cdot \lambda)$ and $x_i \in W_p$. 
  Now, due to the exactness of the translation functors and \cite[Proposition II.7.13]{jan2003}, 
  $T_{\lambda}^{\mu}T(\lambda)$ has a filtration
  \[
    0 = F'_0 \subset \cdots \subset F'_m = T^{\mu}_{\lambda} T(\lambda),
  \]
  where 
   $F'_i/F'_{i-1} = T^{\mu}_{\lambda}V(x_i\cdot \lambda) = V(x_i\cdot \mu)$.
  Likewise,  $\Theta\, T(\lambda)$ has a filtration
  \[
    0 = F''_0 \subset \cdots \subset F''_m = \Theta\, T(\lambda),
  \]
  where $F''_i/F''_{i-1} = T_{\mu}^{\lambda}V(x_i\cdot \mu)$.
  If  $\text{Stab}_{W_p}(\mu) = \{y_1,\dots,y_r\}$, then by \cite[Proposition II.7.13]{jan2003},
  each $T_{\mu}^{\lambda}V(x_i\cdot \mu)$ has a filtration whose layers are $V(y_jx_i\cdot \lambda')$
  for $j=1,\dots,r$. 
  Thus, $\Theta\, T(\lambda')$ has a filtration with layers  $V(y_jx_i\cdot \lambda')$  for $i = 1,\dots,m$ and $j=1,\dots, r$. 
  
  In particular, every layer of the filtration is of the form $V(z\cdot \lambda')$,
  where $z\cdot \mu \leq \mu$. It follows that  $\lambda'$ is a maximal weight with respect to the lattice ordering, and since $\Theta\, T(\lambda)$ is a tilting module, it also follows that
  $T(\lambda') \mid  \Theta\, T(\lambda)$.
\end{proof}

%
\subsection{}
For the rest of this section, assume $p\geq h$ and  $G$ is simple and simply connected.
 Then the ordering $\preceq$ can be directly related to the inclusion ordering of thick tensor ideals. 
\begin{lem}\label{lem:linkage_bound}
 If $\lambda, \mu \in X(T)_+$ satisfy
$C(\lambda) \preceq C(\mu)$, then
\[
\langle T(\lambda) \rangle \supseteq \langle T(\mu) \rangle.
\]
\end{lem}
\begin{proof}
By Proposition~\ref{prop:translation1}, we may assume $\lambda, \mu \in W_p\cdot 0$, since any two
indecomposable tilting modules whose highest weights lie in the lower closure of the same alcove must generate the 
same thick tensor ideal. 
Since $C(\lambda) \preceq C(\mu)$, then by Remark~\ref{rem:super-linkage},
there is a sequence of dominant alcoves 
   \[
   C(\lambda) = C'_0 \uparrow   C'_1 \uparrow \cdots  \uparrow C'_{m_2} = C(\mu),
   \]
   such that  for $i=0,\dots, m_2-1$,
   $C'_{i+1} = s_{\beta_i, n_ip}\cdot C'_i$, where $\beta_i \in \Phi^+$ and 
   $H_{\beta_i, n_ip}$ is an upper wall of $C'_i$. 
   
For $i=0,\dots, m_2-1$, let $\lambda_{i+1} = s_{\beta_i,n_ip}\cdot \lambda_i$, then
\[
\lambda = \lambda_0 \uparrow \lambda_1 \uparrow \cdots \uparrow  \lambda_{m_2} = \mu. 
\] 
Now since $p \geq h$, then it follows from \cite[II.6.3 (1)]{jan2003}
  that for each $i$  there exists an element $\nu_i \in \widecheck{C'_{i}}\cap X(T)_+$  such that 
   $\operatorname{Stab}_{W_p}(\nu_i) = \{1, s_{\beta_i,n_ip}\}$. 
   If  $\Theta_i = T_{\nu_i}^0T_0^{\nu_i}$, then 
 by Proposition~\ref{lem:wall_crossing},  $$T(\lambda_{i+1})\mid \Theta_i T(\lambda_i),$$ and thus
\[
 \langle T(\lambda_i)\rangle \supseteq \langle \Theta_{i}T(\lambda_i)\rangle \supseteq 
 	\langle T(\lambda_{i+1})\rangle,
\]
for $i=0,\dots, m_2-1$. 
Therefore, $ \langle T(\lambda) \rangle \supseteq \langle T(\mu) \rangle.$


\end{proof}

The relationship between thick tensor ideals of tilting modules and their support varieties gives us
a useful corollary. 
\begin{cor}\label{cor:sup-cor}
If $\lambda, \mu \in X(T)_+$ satisfy
$C(\lambda) \preceq C(\mu)$, then
$$V_{G_1}(T(\lambda)) \supseteq V_{G_1}(T(\mu)).$$

\end{cor}
\begin{rem}\label{rem:link-counter}
It would be tempting to hope $\lambda \uparrow \mu$ implies
$V_{G_1}(T(\lambda))\supseteq V_{G_1}(T(\mu))$, but this is not
true in general.  For example, take $G=SL_3(k)$ with $p\geq 3$ 
and let 
$$\lambda + \rho = (p+1)\omega_1 + (p+1)\omega_2,$$
\[
  \mu+\rho = s_{\epsilon_1-\epsilon_2,2p}( \lambda+\rho) 
                =  (3p-1)\omega_1 + 2\omega_2
\]
(see Section~\ref{sec:cell_regions} for notation). 
Then $\lambda \uparrow \mu$, since 
$\mu-\lambda = (p-1)(\epsilon_1-\epsilon_2)$  but 
$\langle \mu + \rho, \epsilon_2-\epsilon_3\rangle < p$, while 
$\langle \lambda + \rho, \epsilon_2-\epsilon_3\rangle > p$. Also, 
$V_{G_1}(T(\lambda)) = \{0\}$, while 
$V_{G_1}(T(\mu))\neq \{0\} $.
\end{rem}

The next lemma relates the support varieties of tilting modules to the support varieties of induced modules. 
\begin{lem}\label{lem:sup-block}
Let $\lambda \in X(T)_+$,
then 
$
V_{G_1}(T(\lambda))\subseteq V_{G_1}(H^0(\mu))
$
for any $\mu \in \widecheck{C(\lambda)}\cap X(T)_+$.
\end{lem}
\begin{proof}
Since $\mu \in \widecheck{C(\lambda)}$, then by Proposition~\ref{prop:translation1}, 
$\langle T(\lambda) \rangle = \langle T(\mu) \rangle$. 
Furthermore,  since $T(\mu)$ has a good filtration whose layers are of the form $H^0(w\cdot \mu)$ with $w \in W_p$,
then by 
 \cite[Proposition 6.2.1]{npv2002}, $V_{G_1}(T(\mu))\subseteq V_{G_1}(H^0(\mu))$. 
 Therefore, 
\[
V_{G_1}(T(\lambda))=V_{G_1}(T(\mu)) \subseteq V_{G_1}(H^0(\mu)).
\]
\end{proof}

We have established the following proposition, which 
 will be a key component in the proof of Theorem~\ref{thm:main}.
\begin{prop}\label{cor:ind-tilting}
Let $\lambda, \mu \in X(T)_+$ satisfy 
$C(\lambda) \preceq C(\mu)$, then for any $\nu \in \widecheck{C(\lambda)}\cap X(T)_+$
\[
 V_{G_1}(H^0(\nu))\supseteq V_{G_1}(T(\lambda))\supseteq V_{G_1}(T(\mu)).
 \]
\end{prop}

\section{Cell regions in type $A_n$}\label{sec:cell_regions}
\subsection{}
For the next two sections, we will assume that $k$ is an algebraically closed field of characteristic $p > 0$ and 
$G=SL_{n+1}(k)$. The roots are given by
\begin{align*}
\Phi  &=\{\epsilon_i-\epsilon_j \mid 1 \leq i,j \leq n+1, \,\, i\neq j\}, \\
\Phi^+ &=\{\epsilon_i-\epsilon_j \mid 1 \leq i<j \leq n+1, \,\, i\neq j\}, \\
\Delta &= \{\epsilon_1-\epsilon_{2}, \dots, \epsilon_{n}-\epsilon_{n+1} \},
\end{align*}
where $\epsilon_1,\dots, \epsilon_{n+1}$ is the standard basis for $\mathbb{E} \cong \mathbb{R}^{n+1}$. 
 The corresponding fundamental weights are denoted by $\omega_1,\dots,\omega_n \in X(T)_+$. 
Since $\Phi$ is simply-laced, it can be normalized so that every root has length $\sqrt{2}$. Then $\alpha^{\vee} = \alpha$ for all 
$\alpha \in \Phi$. The Weyl group $W = \Sigma_{n+1}$ is the group of permutations of $\{1,\dots, n+1\}$, where for any 
$w \in \Sigma_{n+1}$ 
 and $\epsilon_i-\epsilon_j \in \Phi$, define
\[
w(\epsilon_i-\epsilon_j) = \epsilon_{w(i)}-\epsilon_{w(j)}.
\]
In particular, the reflections $s_{\epsilon_i-\epsilon_j} = (i,j) \in \Sigma_{n+1}$ (in cycle notation). 

Let $\mathcal{P}$ denote the set of all partitions of $n+1$.
For each partition $\pi = (p_1,p_2,\dots, p_r)  \in \mathcal{P}$ with $p_1\geq p_2 \geq \cdots \geq p_r \geq 1$, let  
 $x_{\pi} \in \mathcal{N}$ denote the nilpotent matrix which is a direct sum of Jordan blocks of sizes 
 $p_1,\dots, p_r$ and let $\mathcal{O}_{\pi} = G\cdot x_{\pi}$.
The assignment $\pi \leftrightarrow \mathcal{O}_{\pi}$ gives a bijection between $\mathcal{P}$ and the set of $G$-orbits in $\mathcal{N}$ (\cite{cm}). 

 \begin{defn}\label{defn:dom-order}
The set $\mathcal{P}$ 
is equipped with the 
\emph{dominance ordering}, where if
$\pi = (p_1,p_2,\dots, p_{r})$ and $\sigma = (q_1,q_2,\dots,q_{s})$, then 
$\pi \leq \sigma$ if and only if for $k=1,\dots, n+1$, 
\[
p_1 + \cdots + p_k \leq q_1 + \cdots + q_k,
\] 
where $p_k=0$ (resp. $q_k=0$) if $k > r$ (resp. $k > s$). 
Moreover, $\mathcal{P}$ has an order reversing \emph{transposition} operation, denoted 
$\pi \mapsto \pi^t$. 
\end{defn}
\begin{rem} 
If $\pi$, $\sigma \in \mathcal{P}$, then $\overline{\mathcal{O}_{\pi}} \subseteq \overline{\mathcal{O}_{\sigma}}$ if and only if 
 $\pi \leq \sigma$ (cf. \cite{cm}). 
 \end{rem} 
Furthermore, $\mathcal{P}$ also has a \emph{supremum} (or \emph{least upper bound})
operation with respect to $\leq$.
\begin{defn}\label{defn:sup}
Let $\{\pi_1,\dots, \pi_t\}$ be any subset of $\mathcal{P}$, then there exists a \emph{least upper bound}
$\pi =  \text{sup}\{\pi_1,\dots,\pi_t\}.$
 To define $\pi$, set $\pi_i = (p_{i,1},\dots, p_{i,r_i})$ for $i=1, \dots, t$, then 
$\pi = (p_1,\dots, p_{n+1})$, where $p_1 = \text{max}\{p_{1,1},\dots, p_{t,1}\}$ and for $i\geq 2$,
\[
p_1 + \cdots + p_i = \text{max}\{p_{1,1} + \cdots + p_{1,i}, \dots, p_{t,1} + \cdots + p_{t,i}\}.
\]
It follows from Definition~\ref{defn:dom-order} that $\pi$ satisfies the least upper bound property. 
%
\end{defn}
\begin{exmp}
Suppose $n+1 = 6$, $\pi_1 = (3,3)$ and $\pi_2 = (4,1,1)$, then 
$
\text{sup}\{\pi_1,\pi_2\} = (4,2).
$
\end{exmp}

From Definition~\ref{defn:good-prime}, recall  that every prime $p$ is good for $\Phi$.  As $d_{\alpha} = 1$ for all $\alpha$, we have
 \[
 \Phi_{\lambda,p} = \{ \alpha \in \Phi \mid \langle \lambda + \rho, \alpha \rangle \in p\mathbb{Z}\}
 \]
for any $\lambda \in \mathbb{E}$ (see Definition~\ref{defn:stab-ss}).


\begin{defn}\label{defn:root-partition}
 Any subsystem $\Phi' \subseteq \Phi$ is conjugate to one of type
$A_{p_1-1}\times A_{p_2-1}\times\cdots A_{p_{r}-1}$,
where $p_1\geq p_2\geq \cdots \geq p_r \geq 1$ and $A_0$ denotes the type of the empty root system.
Thus, it can be associated to a partition which is given by
\[
\pi(\Phi')= (p_1,p_2,\dots, p_r) \in \mathcal{P}.
\]
More explicitly,  this partition is obtained by choosing a basis $\Delta(\Phi')$ for $\Phi'$ and 
decomposing $\Delta(\Phi') = \Delta_1 \sqcup \cdots \sqcup \Delta_r$, so that for some $w \in W$, 
\[
w(\Delta(\Phi')) = I_1 \sqcup I_2\sqcup \cdots \sqcup I_r \subseteq \Delta,
\]
where $w(\Delta_k) = I_k$,  $\Phi_{I_k}$ is of type  $A_{p_k-1}$  and
$p_k = |\Delta_k| +1$  for $k=1,\dots, r$ (if $p_k =1$, then $\Delta_k = I_k = \emptyset$ and
$\Phi_{I_k} =\emptyset$).  
\end{defn}

 This allows us to associate a partition to each 
$\lambda \in \mathbb{E}$ (or equivalently to each $F \in \mathcal{F}$). 
\begin{defn}\label{def:weyl-part}
 For any $\lambda \in \mathbb{E}$,
  let $d(\lambda) = \pi(\Phi_{\lambda,p})$. 
 \end{defn}
 If $\lambda \in X(T)_+$, then by \cite[Theorem 6.2.1]{npv2002},
  $V_{G_1}(H^0(\lambda)) = \overline{\mathcal{O}_{d(\lambda)^t}}$.

\subsection{}
Now we shall give a description of the cell regions in type $A$ which is inspired by the treatment given 
by Shi in \cite{shi} and paraphrased by Cooper in 
\cite[4.1]{coop2010}.
 \begin{defn}\label{defn:good-subsystem}
A subset $\Psi \subseteq \Phi^+$ is said to be a \emph{positive subroot system} 
if $\Psi=  w(\Phi_I^+)$ for some $I \subseteq \Delta$ and $w \in W$. The set 
$\Delta(\Psi) = w(I)$ will be referred to as a \emph{basis} of $\Psi$ because $w(I)$ is a basis
of the subroot system $w(\Phi_I) \subseteq \Phi$. We can define $\pi(\Psi) = \pi(\Phi_I)$ to be the 
partition associated to $\Psi$. 
\end{defn}

A basis of a positive subroot system is also a basis for an actual subroot system of $\Phi$.
Thus,
a subset $\Delta' \subseteq \Phi^+$ is a basis for a positive subroot system only if 
for any two distinct elements $\alpha,\beta \in \Delta'$ with
 $\alpha = \epsilon_i - \epsilon_j$ and  $\beta=\epsilon_k - \epsilon_l$, one has
 $m_{\alpha,\beta}=\langle \alpha, \beta \rangle \in \{0, -1\}$. The $m_{\alpha,\beta}=0$ case
 occurs precisely when the indices $i,j,k,l$ are all distinct, and the $m_{\alpha,\beta}=-1$ case occurs
 precisely when either $i=l$ or $j=k$. 
 Conversely, suppose $\Delta'\subseteq \Phi^+$ is a basis for a subroot system of $\Phi$, then 
$\Delta' = \Delta_1\sqcup \cdots \sqcup \Delta_r$ with
  $$\Delta_k = \{ \epsilon_{i_{k,1}}-\epsilon_{i_{k,2}}, \epsilon_{i_{k,2}}-\epsilon_{i_{k,3}}, \dots
   \epsilon_{i_{k,p_k-1}}-\epsilon_{i_{k,p_k}}\},$$ 
   where 
   $i_{k,1} < i_{k,2} < \cdots < i_{k,p_k}$ and $p_1 \geq p_2 \geq \cdots \geq p_r \geq  2$.
   Let $w \in W = \Sigma_{n+1}$, where
   \[
   w(i_{k,j}) = p_1 + p_2 + \cdots + p_{k} + j
   \]
   for $k=1,\dots, r$ and $j=1,\dots, p_k$, and 
   \[
   w(i) = i
   \]
   if $i \neq i_{k,j}$ for some $k,j$. Now if $I = I_1\sqcup I_2\sqcup  \cdots \sqcup I_r$, where 
   \[
   I_k = \{\epsilon_{p_1 + \cdots + p_{k-1} +j} - \epsilon_{(p_1 + \cdots + p_{k-1} +j)+1} \}_{j=1,\dots, p_k-1}
   \]
  for $k=1,\dots, r$,   then $\Delta' = w^{-1}(I)$. Thus, $\Delta'$ is a basis for the positive subroot system 
   $\Psi = w^{-1}(\Phi_{I}^+)$ and $\Delta' = \Delta(\Psi)$.  Moreover, the partition associated to $\Psi$
   is given by
   \[
   \pi(\Psi) = (p_1,\dots,p_r,1,\dots,1).
   \]

 For example, if $n+1 =7$, then the subset
 $\{\epsilon_1-\epsilon_6, \epsilon_2-\epsilon_4, \epsilon_4-\epsilon_7\}\subseteq \Phi^+$ is a basis 
 for a positive subroot system. However, the subset 
 $\{\epsilon_3-\epsilon_6,\epsilon_3-\epsilon_5\} \subseteq \Phi^+$ 
 is not a basis for a positive subroot system. 

The following partition can be used to describe both the weight cells and the Lusztig bijection in type 
$A_n$
(cf. Section~\ref{sec:weight-cells}). It was originally defined by Shi (cf. \cite{shi}). However, the following formulation
of the definition was
 given by Cooper in \cite[4.1]{coop2010}. 
\begin{defn}\label{def:hum-part}
For each $\lambda \in \mathcal{C}$, set 
\[
\Gamma_{\lambda}= \{\alpha  \in \Phi^+ 
		\mid \langle \lambda + \rho,  \alpha \rangle \geq p \} 
\]
and define 
\[
s(\lambda) = \text{sup}\{ \pi(\Psi) \mid \text{$\Psi \subseteq \Gamma_{\lambda}$ is a positive
		subroot system}\}.
\]
 \end{defn}
\begin{rem}\label{rem:weight-cell}
 Under the correspondence between $\mathcal{P}$ and the set of nilpotent orbits, 
 each partition $\pi \in \mathcal{P}$ defines a \emph{weight cell} $c_{\pi} \subseteq \mathcal{C}$, consisting
of all $\lambda \in \mathcal{C}$ satisfying $s(\lambda)^t = \pi$. The bijection
$c_{\pi} \leftrightarrow \mathcal{O}_{\pi}$ establishes the Lusztig bijection between weight cells 
and nilpotent orbits. 
\end{rem}

 At this point, all of the notation required to understand the 
 statement of Theorem ~\ref{thm:main} has been introduced. 
In fact, a more general conjecture which places no assumption on $p$ was formulated, and then verified for $p=2$ by
Cooper in \cite{coop2010}. 
 \begin{conj}[Cooper]\label{conj:hump}
 For any $\lambda \in X(T)_+$, 
 $V_{G_1}(T(\lambda)) = \overline{\mathcal{O}_{s(\lambda)^t}}$.
 \end{conj}

Conjecture~\ref{conj:hump} is equivalent to the statement
that for  $\pi \in \mathcal{P}$ and $\lambda \in c_{\pi}\cap X(T)_+$, 
$V_{G_1}(T(\lambda))=\overline{\mathcal{O}_{\pi}}.$

\subsection{Good positive subroot systems}\label{sec:weight-cell}

According to Definition~\ref{def:hum-part}, $s(\lambda)$ is calculated by taking the supremum 
over all partitions of the form $\pi(\Psi)$, where $\Psi \subseteq \Gamma_{\lambda}$ is a positive 
subroot system of $\Gamma_{\lambda}$.
 However, it will soon be shown that $s(\lambda)$ can also be calculated by taking the supremum of a smaller subset
 of partitions. Namely, the set of partitions of the form $\pi(\Psi)$, where
 $\Psi \subseteq \Gamma_{\lambda}$ is a \emph{good} positive subroot system (see the following definition). 
%

 \begin{defn}
 A positive subroot system $\Psi \subseteq \Phi^+$ is called \emph{good}, if 
 there are no two elements $\alpha, \beta \in \Delta(\Psi)$ satisfying 
 $\alpha < \beta$. 
 \end{defn}
 \begin{rem}
 If $\Psi \neq \emptyset$, then it can be verified that 
 $\Psi$ is good if and only if 
 \[
 \Delta(\Psi) = \{\epsilon_{i_1}-\epsilon_{j_1}, \epsilon_{i_2}-\epsilon_{j_2},\dots, \epsilon_{i_r}-\epsilon_{j_r}\},
 \]
 where $i_1<i_2<\cdots <i_r$ and $j_1<j_2<\cdots <j_r$. 
 \end{rem}
 
 For example, if $n+1 = 5$, then the positive subroot system with basis 
 $\{\epsilon_1-\epsilon_4, \epsilon_2-\epsilon_3\}$ is not good since 
 $\epsilon_2-\epsilon_3 < \epsilon_1 - \epsilon_4$. On the other hand,
 $\{\epsilon_1-\epsilon_4, \epsilon_2 - \epsilon_5\}$ is the basis for 
 a good positive subroot system. 

For any positive subroot system $\Psi \subseteq \Phi^+$, let
\[
\Gamma_{\Psi} = \{ \alpha \in \Phi^+ \mid \alpha \geq \beta \text{ for some $\beta \in \Psi$}\}.
\]
If $\Psi \subseteq \Gamma_{\lambda}$ is a positive subroot system, then 
 $\Gamma_{\Psi} \subseteq \Gamma_{\lambda}$ because for any $\alpha \in \Gamma_{\Psi}$, there exists 
 $\beta \in \Psi \subseteq \Gamma_{\lambda}$ satisfying $\alpha \geq \beta$. Thus, since $\lambda \in \mathcal{C}$, 
 \[
 \langle \lambda + \rho, \alpha \rangle \geq \langle \lambda + \rho, \beta \rangle \geq p.
 \] 
 Moreover, for any good positive subroot system
  $\Psi' \subseteq \Gamma_{\Psi}$,  $\pi(\Psi') \leq s(\lambda)$. 

\begin{lem}\label{lem:goodsubsystem}
Let $\Psi \subseteq \Phi^+$ be a positive subroot system. Then there exist
good positive subroot systems $\Psi_1,\dots, \Psi_t \subseteq \Gamma_{\Psi}$ such that 
$\pi(\Psi) \leq \text{sup}\{\pi(\Psi_1),\dots, \pi(\Psi_t)\}$. 
\end{lem}
\begin{proof}
For each positive subroot system $\Psi\subseteq \Phi^+$, let
\[
M_{\Psi}= \{ \{\alpha,\beta \} \subseteq \Delta(\Psi) \mid \alpha < \beta  \}
\] and
$
m_{\Psi} = |M_{\Psi}|.
$
We shall perform induction on $m=m_{\Psi}\geq 0$. The base case, $m=0$, follows from the
fact that  
$m_{\Psi} = 0$ if and only if $\Psi$ is good. 
For the inductive step,  suppose $m \geq 1$ is arbitrary and that for any 
positive subroot system $\Psi \subseteq \Phi^+$ satisfying $m_{\Psi} < m$, 
there exist
good positive subroot systems $\Psi_1,\dots, \Psi_t \subseteq \Gamma_{\Psi}$ such that 
$\pi(\Psi) \leq \text{sup}\{\pi(\Psi_1),\dots, \pi(\Psi_t)\}$. 

Now suppose $\Psi \subseteq \Phi^+$ satisfies $m_{\Psi} = m$. Let
$\Delta(\Psi) = \Delta_1\sqcup \cdots \sqcup \Delta_r$ and $\pi(\Psi) = (p_1,\dots, p_r)$ be as in 
Definitions~\ref{defn:root-partition} and \ref{defn:good-subsystem}.
 For each $k$ satisfying  $\Delta_{k} \neq \emptyset$ (i.e., $p_k \geq 2$), write
  $$\Delta_k = \{ \epsilon_{i_{k,1}}-\epsilon_{i_{k,2}}, \epsilon_{i_{k,2}}-\epsilon_{i_{k,3}}, \dots
   \epsilon_{i_{k,p_k-1}}-\epsilon_{i_{k,p_k}}\},$$ 
   where 
   $i_{k,1} < i_{k,2} < \cdots < i_{k,p_k}$. 
   Since $m_{\Psi} \geq 1$, there exists $\alpha_1 = \epsilon_{i_{t_1,s_1}}-\epsilon_{i_{t_1,s_1+1}} \in \Delta_{t_1}$ 
   and  $\alpha_2 = \epsilon_{i_{t_2,s_2}} -\epsilon_{i_{t_2,s_2+1}} \in \Delta_{t_2}$ such that
    $\alpha_1 > \alpha_2$. Our goal is to construct a new positive subroot system $\Psi''\subseteq \Phi^+$ 
    by replacing the two 
    ``bad'' roots $\alpha_1,\alpha_2$ 
   with two non-comparable roots $\beta_1, \beta_2$ in such a way so that $\Psi' \subseteq \Gamma_{\Psi}$ and 
   $m_{\Psi'} < m_{\Psi}$.
   
   To do this, let $\Psi'\subseteq \Gamma_{\Psi}$ denote the subroot system whose basis, $\Delta(\Psi')$,
   is obtained by taking $\Delta(\Psi)$ and replacing the roots $\alpha_1, \alpha_2$ with
   $\beta_1= \epsilon_{i_{t_1,s_1}}-\epsilon_{i_{t_2, s_2+1}}$ and 
   $\beta_2 = \epsilon_{i_{t_2,s_2}}-\epsilon_{i_{t_1,s_1+1}}$ respectively. 
   The inclusion
   $\Psi'\subseteq \Gamma_{\Psi}$
   holds because  $\beta_k \geq \alpha_2$ for $k=1,2$ (which implies  $\beta_1, \beta_2 \in \Gamma_{\Psi}$). 
   There is a decomposition,
     $\Delta(\Psi') = \Delta'_1\sqcup \cdots \sqcup \Delta'_r,
     $
     where  
     \begin{align*}
     \Delta'_{t_1} &= \{ \epsilon_{i_{t_1,1}}-\epsilon_{i_{t_1,2}}, \dots,  \epsilon_{i_{t_1,s_1}}-\epsilon_{i_{t_2,s_2+1}}, \dots,
   \epsilon_{i_{t_2,p_{t_2}-1}}-\epsilon_{i_{t_2,p_{t_2}}}\}, \\
    \Delta'_{t_2} &= \{ \epsilon_{i_{t_2,1}}-\epsilon_{i_{t_2,2}}, \dots,  \epsilon_{i_{t_2,s_2}}-\epsilon_{i_{t_1,s_1+1}}, \dots,
   \epsilon_{i_{t_1,p_{t_1}-1}}-\epsilon_{i_{t_1,p_{t_1}}}\}
     \end{align*}
     and $\Delta'_k = \Delta_k$ for $k \not\in \{t_1,t_2\}$.
     For each $k$, let $p'_k = |\Delta'_k|+1$.
      Then $p'_k = p_k$ for $k \not\in \{t_1,t_2\}$ and $p'_{t_1} + p'_{t_2} = p_{t_1} + p_{t_2}$
       (however, we cannot assume that 
     $p'_k \geq p'_j$ whenever $k < j$). 
     In any case,
     \[
     \pi(\Psi') = (p'_{\tau(1)},p'_{\tau(2)}, \dots, p'_{\tau(r)}),
     \]
     where $p'_{\tau(1)}\geq p'_{\tau(2)}\geq \cdots \geq p'_{\tau(r)}$ and
     $\tau$ is a permutation of $\{1,2,\dots, r\}$.   Also,
     \[
     p'_1 + \cdots + p'_i \leq p'_{\tau(1)} + \cdots + p'_{\tau(i)}
     \]
     for $i=1,\dots,r$. 


   Furthermore, $m_{\Psi'} < m_{\Psi}$.  To see why this is true, 
  begin by observing  that if 
    $$\{\alpha,\beta\} \subseteq 
     \Delta(\Psi')\backslash \{\beta_1,\beta_2\}=\Delta(\Psi)\backslash \{\alpha_1,\alpha_2\},$$ 
   then $\{ \alpha,\beta\} \in M_{\Psi'}$ if and only if $\{\alpha,\beta\} \in M_{\Psi}$. By definition,
   $\{\beta_1,\beta_2\} \not\in M_{\Psi'}$.
   Thus, it will be sufficient to show that 
   the number of subsets of the form $\{\alpha,\beta_{k}\} \in M_{\Psi'}$ is no greater than the number
   of subsets of the form $\{\alpha,\alpha_k\} \in M_{\Psi}$,
   where $\alpha \in  \Delta(\Psi')\backslash \{\beta_1,\beta_2\}$. 
   
    First suppose $\alpha \in \Delta(\Psi')$ satisfies  $\{\alpha,\beta_k\} \in M_{\Psi'}$ for both $k=1,2$. 
   If $\alpha > \beta_k$ for $k=1,2$, then $\alpha > \alpha_1 > \alpha_2$, and hence
   $\{\alpha,\alpha_k\} \in M_{\Psi}$ for both $k$. Similarly, if $\alpha < \beta_k$
   for both $k$, then $\{\alpha,\alpha_k\} \in M_{\Psi}$ for both $k$. If $\alpha > \beta_1$ and 
   $\alpha < \beta_2$, then $\alpha > \alpha_2$ and $\alpha < \alpha_1$ (since $\beta_1> \alpha_2$
   and $\beta_2 < \alpha_1$), and hence $\{\alpha,\alpha_k\} \in M_{\Psi}$ for $k=1,2$. Likewise, if 
   $\alpha < \beta_1$
   and $\alpha > \beta_2$, then $\{\alpha,\alpha_k\} \in M_{\Psi}$ for both $k$. 
   Suppose now that $\{\alpha,\beta_k\} \in M_{\Psi'}$ only for a single $k$. In this case, if
    $\alpha < \beta_k$, then $\alpha < \alpha_1$ since $\beta_k < \alpha_1$, and thus
    $\{\alpha,\alpha_1\} \in M_{\Psi}$. Similarly, if $\alpha > \beta_k$, then $\alpha > \alpha_2$ and 
     $\{\alpha,\alpha_2\} \in M_{\Psi}$. 
    So the number of pairs of the form $\{\alpha,\beta_k\} \in M_{\Psi'}$ is
    no greater than the number of pairs of the form $\{\alpha,\alpha_k\} \in M_{\Psi}$.  
    Therefore, $m_{\Psi'} \leq m_{\Psi}-1$. 
    
     
     For simplicity, assume $t_1 < t_2$ (the exact same argument will also work when 
     $t_1 > t_2$).  
     Define $\Psi''$ to be the positive subroot system with basis
     $\Delta(\Psi'')= \Delta(\Psi)\backslash \Delta_{t_2}$. Then
     $\Psi'' \subseteq \Gamma_{\Psi}$
     and  $m_{\Psi''} < m_{\Psi}$, since $M_{\Psi''} \subseteq M_{\Psi}\backslash \{\{\alpha_1,\alpha_2\}\}$.
     The partition associated to $\Psi''$ is given by
      \begin{equation*}\label{eqn:pi2}
      \pi(\Psi'') = (p_1, \dots, p_{t_2-1}, p_{t_2+1}, \dots,p_r,1,\dots,1).
      \end{equation*}
     
     Let $\pi = \text{sup}\{\pi(\Psi'),\pi(\Psi'')\}$, and write $\pi = (q_1,\dots, q_{n+1})$ with 
     $q_1\geq q_2 \geq \cdots \geq q_{n+1}\geq 0$. 
     By Definition~\ref{defn:sup}, 
     $\pi(\Psi'') \leq \pi$, and hence
     \[
      p_1 + \cdots + p_i \leq q_1 + \cdots + q_i
     \]
     for $i=1,\dots, t_2-1$. Moreover, since $p'_{t_1}+p'_{t_2} = p_{t_1} + p_{t_2}$, then 
     \[
     p_1 + \cdots + p_i = p'_1 + \cdots + p'_i \leq p'_{\tau(1)} + \cdots + p'_{\tau(i)} \leq q_1 + \cdots + q_i
     \]
     for $i\geq t_2$. Thus, $\pi(\Psi) \leq \pi$. Now by the inductive hypothesis, since $m_{\Psi'} < m$ and 
     $m_{\Psi''} < m$, there exist good positive subroot systems $\Psi'_1,\dots, \Psi'_{t'} \subseteq \Gamma_{\Psi'}$ and
     $\Psi''_1,\dots, \Psi''_{t''} \subseteq \Gamma_{\Psi''}$ such that 
     $\pi(\Psi') \leq \text{sup}\{\pi(\Psi'_1),\dots, \pi(\Psi'_{t'})\}$ and 
     $\pi(\Psi'') \leq \text{sup}\{\pi(\Psi''_1),\dots, \pi(\Psi''_{t''})\}$. It follows that 
     \[
     \pi(\Psi) \leq \pi \leq \text{sup}\{\pi(\Psi'_1),\dots,\pi(\Psi'_{t'}),\pi(\Psi''_{1}),\dots, \pi(\Psi''_{t''})\},
     \]
     where $\Psi'_1,\dots,\Psi'_{t'},\Psi''_1,\dots, \Psi''_{t''} \subseteq \Gamma_{\Psi}$ are good positive 
     subroot systems. 

\end{proof}
For the sake of clarity, a nontrivial example demonstrating the algorithm which was used in the above lemma has 
been included. 
\begin{exmp}
Suppose that $n+1=6$ and let $\Delta(\Psi) = \Delta_1\sqcup \Delta_2$, where
\begin{align*}
\Delta_1 &= \{\epsilon_1-\epsilon_3, \epsilon_3-\epsilon_4, \epsilon_4-\epsilon_6\} \\
\Delta_2 &= \{\epsilon_2-\epsilon_5\}.
\end{align*}
Then $m_{\Psi} = 1 = |\{ \{\epsilon_3-\epsilon_4, \epsilon_2-\epsilon_5\}\}|$ and $\pi(\Psi) = (4,2)$. 
Let $\alpha_1 = \epsilon_2-\epsilon_5$ and $\alpha_2 = \epsilon_3-\epsilon_4$ be as in the above proof, then 
$\beta_1 = \epsilon_2-\epsilon_4$ and $\beta_2 = \epsilon_3-\epsilon_5$. Thus, 
$\Delta(\Psi_1) = \Delta'_1\sqcup \Delta'_2$, where
\begin{align*}
\Delta'_1 &= \{\epsilon_1-\epsilon_3, \epsilon_3-\epsilon_5\}\\ 
\Delta'_2 &=  \{\epsilon_2-\epsilon_4,\epsilon_4-\epsilon_6\}.
\end{align*}
Then $\Psi_1$ is good, $\Gamma_{\Psi_1}\subseteq \Gamma_{\Psi}$ and 
$\pi(\Psi_1) = (3,3) \leq \pi(\Psi)$.  Following the algorithm given in the preceding proof, we can also obtain
the good subroot system $\Psi_2$ with basis
\begin{align*}
\Delta(\Psi_2) = \{\epsilon_1-\epsilon_3, \epsilon_3-\epsilon_4, \epsilon_4-\epsilon_6\},
\end{align*}
by removing $\epsilon_2-\epsilon_5$ from $\Delta(\Psi)$. Then 
$\Gamma_{\Psi_2} \subseteq \Gamma_{\Psi}$ and $\pi(\Psi_2) = (4,1,1)$. Finally, observe that 
$\pi(\Psi) =\text{sup}\{\pi(\Psi_1),\pi(\Psi_2)\}$.

\end{exmp}

The following proposition gives an equivalent characterization of
the partitions $s(\lambda) \in \mathcal{P}$ for $\lambda \in \mathcal{C}$. 
\begin{prop}\label{cor:good-part}
For each $\lambda \in \mathcal{C}$, 
\[
s(\lambda) = \operatorname{sup}\{ \pi(\Psi) \mid \text{$\Psi \subseteq \Gamma_{\lambda}$ is a good positive
		subroot system}\}.
\] 
\end{prop}
\begin{proof}
Let $\pi =  \text{sup}\{ \pi(\Psi) \mid \text{$\Psi \subseteq \Gamma_{\lambda}$ is a good positive
		subroot system}\}.$
It follows immediately that $\pi \leq s(\lambda)$. Conversely, by Lemma~\ref{lem:goodsubsystem}, 
for each positive subroot system $\Psi \subseteq \Gamma_{\lambda}$, there exist
good positive subroot systems 
$
\Psi_1,\dots, \Psi_t \subseteq \Gamma_{\Psi} \subseteq \Gamma_{\lambda}
$
such that 
\[
\pi \geq \text{sup}\{\Psi_1,\dots, \Psi_t\} \geq \pi(\Psi).
\]
Thus, $\pi \geq \pi(\Psi)$ for any positive subroot system $\Psi \subseteq \Gamma_{\lambda}$, and hence 
$\pi \geq s(\lambda)$. 
\end{proof}

\section{The upper bound}\label{sec:upper-bound}
\subsection{}
 In this section, we will 
  prove the upper bound portion of Theorem~\ref{thm:main}.
   A key tool will be the following lemma, which illustrates the
    importance of good positive subroot systems. 
\begin{lem}\label{lem:gamma-subroot}
If $\lambda \in \mathcal{C}$ and
$\Psi \subseteq \Gamma_{\lambda}\subseteq \Phi^+$ is a good
positive subroot system,  then there exists an element
$\mu \in \mathcal{C}$ such that  $\Phi_{\mu,p}^+ \supseteq \Psi$ and 
 $C(\mu) \preceq C(\lambda)$. 
\end{lem}
\begin{proof}
Let $\lambda \in \mathcal{C}$ and let $\Psi \subseteq \Gamma_{\lambda}$ be a good 
positive subroot system. 
If $\Psi = \emptyset$, let  
$\mu + \rho = (1/n,\dots, 1/n)$ be given in fundamental basis coordinates. Then $C(\mu)$
is the bottom alcove, and hence $C(\mu) \preceq C(\lambda)$ since $C(\lambda) \in \mathcal{A}^+$. 
Furthermore, $\Phi_{\mu,p}^+ \supseteq \Psi = \emptyset$. 

For the rest of the proof, we shall assume that $\Psi \neq \emptyset$. In this case, 
$\Delta(\Psi) = \{\alpha_1,\dots, \alpha_r \}$, where $\alpha_k = \epsilon_{i_k}-\epsilon_{j_k}$,
$i_1<i_2<\cdots < i_r$ and  $j_1 < j_2 < \cdots < j_r$.  
 By performing induction on the \emph{rank} $r = |\Delta(\Psi)| \geq 1$, it will be shown that there exists 
$\mu \in \mathcal{C}$ satisfying
 \begin{enumerate} 
\item $\langle \mu + \rho, \alpha_i \rangle =p$  for any $\alpha_i \in \Delta(\Psi)$, and thus $\Phi_{\mu,p}^+ \supseteq \Psi$, \\
\item $\langle \mu + \rho, \epsilon_1-\epsilon_{j_1-1} \rangle < p$  and 
    $\langle \mu + \rho, \epsilon_{i_r+1}-\epsilon_{n+1}\rangle < p$, \\
\item $C(\mu)\preceq C(\lambda)$.
\end{enumerate} 
For each $\epsilon_i-\epsilon_j \in \Phi^+$, let $n_{ij}\geq 1$ be the unique integer satisfying
\begin{align*}
(n_{ij}-1)p \leq \langle \lambda + \rho, \epsilon_i - \epsilon_j \rangle < n_{ij}p.
\end{align*}
It is useful to remark that Lemma~\ref{lem:ss-implies-s} implies that (3) will follow if
\begin{equation}\label{eqn:reduced-problem}
\langle \mu + \rho, \epsilon_i - \epsilon_j \rangle <n_{ij}p
\end{equation}
is satisfied for all $\epsilon_i - \epsilon_j \in \Phi^+$.

For the base case, suppose $r=1$,  then $\Delta(\Psi) = \{\epsilon_{i_1} - \epsilon_{j_1}\}$. 
Now let $\mu \in \mathcal{C}$ be given by 
\[
\mu + \rho = (a_1,\dots,a_n),
\]
where $a_k=1/n$ for $k\geq i_1+1$, $a_{i_1}  = p-\frac{j_1-i_1-1}{n}$ and, if $i_1 > 1$, 
 $a_k = a$ for $k\leq i_1-1$ with 
 $
0< a < \frac{1}{(i_1-1)n}.
 $
Then $\mu$ automatically satisfies (1). Moreover, 
\begin{align*}
\langle \mu + \rho, \epsilon_1-\epsilon_{j_1-1} \rangle &= \langle \mu + \rho, \epsilon_1-\epsilon_{i_1} \rangle + 
		\langle \mu + \rho, \epsilon_{i_1}-\epsilon_{i_1+1} \rangle + 
			\langle \mu + \rho, \epsilon_{i_1+1}-\epsilon_{j_1-1} \rangle \\
			&= (i_1-1)a  +  \left( p-\frac{j_1-i_1-1}{n} \right) + \frac{j_1-i_1 -2}{n} \\
			&= (i_1-1)a + p - \frac{1}{n} < p
\end{align*}
and 
\[
\langle \mu + \rho, \epsilon_{i_1+1}-\epsilon_{n+1} \rangle = \frac{n-i_1-1}{n} < p,
\]
so (2) is also satisfied. To show (3), first observe that since (2) holds, then  if 
$\epsilon_i -\epsilon_j \ngeq \epsilon_{i_1}-\epsilon_{j_1}$, 
\[
\langle \mu + \rho, \epsilon_i - \epsilon_j \rangle < p \leq n_{ij}p. 
\]
On the other hand, since $\langle \mu + \rho, \epsilon_1-\epsilon_{n+1}\rangle < p+1 <2p$, then 
if $\epsilon_i - \epsilon_j \geq \epsilon_{i_1}-\epsilon_{j_1}$, 
\[
\langle \mu + \rho, \epsilon_i - \epsilon_j \rangle < 2p \leq n_{ij}p,
\]
where $n_{ij}\geq 2$ because 
\[
p\leq  \langle \lambda + \rho, \epsilon_{i_1} - \epsilon_{j_1} \rangle 
\leq \langle \lambda + \rho, \epsilon_i - \epsilon_j \rangle,
\]
and hence $2\leq n_{i_1j_1} \leq n_{ij}$. 

For the inductive step, suppose $r \geq 2$ and that for any good positive subroot system 
$\Psi' \subseteq \Gamma_{\lambda}$ with $|\Delta(\Psi')| < r$, there exists $\nu \in \mathcal{C}$ satisfying
conditions (1), (2) and (3). 
Now let $\Psi \subseteq \Gamma_{\lambda}$ be a good positive subroot system with 
$\Delta(\Psi) = \{\alpha_1,\dots, \alpha_r \}$, where $\alpha_k = \epsilon_{i_k}-\epsilon_{j_k}$,
$i_1<i_2<\cdots < i_r$ and  $j_1 < j_2 < \cdots < j_r$. 
 Our goal is to find an element $\mu \in \mathcal{C}$ corresponding to $\Psi$ which satisfies all three conditions. 
Suppose
\[
\mu + \rho = (a_1, \dots, a_n),
\]
using fundamental basis coordinates. 
Thus, to determine $\mu$, it suffices to determine the appropriate coordinates $a_1,\dots, a_n$. 

Let $\Psi' \subset \Psi$ be the good positive subroot system with basis 
$\Delta(\Psi')= \{\alpha_2,\dots, \alpha_r \}$, then by the inductive hypothesis, there exists
$\nu \in \mathcal{C}$ corresponding to $\Psi'$ which is given by 
$$\nu+\rho = (b_1,\dots,b_n),$$
and satisfies conditions (1), (2) and (3). 
%
We begin by choosing the coordinates $a_{i_1},\dots, a_{n}$.
Set $a_k = b_k$ for all $k\geq i_1+1$, and set
\begin{equation*}\label{eqn:condition-0}
a_{i_1} = p-(b_{i_1+1}+\cdots +b_{j_1-1}).
\end{equation*}
Observe that
$b_1 + \cdots + b_{j_1-1} < p$ because $\nu$ satisfies (2) and $j_1 < j_2$.
It follows that $a_{i_1} > 0$. 
Furthermore, 
$$
\langle \mu + \rho, \epsilon_{i_1}-\epsilon_{j_1}\rangle = 
	p- (b_{i_1+1}+\cdots +b_{j_1-1}) + b_{i_1+1} + \cdots + b_{j_1-1}  = p
$$
and for each $k\geq 2$,
\[
\langle \mu + \rho, \epsilon_{i_k}-\epsilon_{j_k}\rangle = 
\langle \nu + \rho, \epsilon_{i_k}-\epsilon_{j_k}\rangle = p.
\]
Moreover,  
\begin{equation*} 
\langle \mu + \rho, \epsilon_{i_r+1} - \epsilon_{n+1}\rangle = \langle \nu + \rho, \epsilon_{i_r+1} - \epsilon_{n+1}\rangle  < p. 
\end{equation*}
Thus, if $a_{i_1},\dots, a_n$ are chosen in this way, then $\mu$ satisfies (1)
 for any choice of positive real numbers $a_1,\dots, a_{i_1-1}$. 
 Additionally,  if $i_1 = 1$, then $\mu$ already satisfies (2) because
\[
\langle \mu + \rho, \epsilon_1-\epsilon_{j_1-1} \rangle = p - b_{j_1-1}  < p. 
\]
Also, if $\epsilon_i - \epsilon_j \in \Phi^+$ and $i \geq i_1+1 =2$, then
\[
\langle \mu + \rho, \epsilon_i-\epsilon_{j} \rangle = \langle \nu + \rho, \epsilon_i-\epsilon_{j} \rangle
 < n_{ij}p,
\]
 If $i =1$ and $j\leq j_1-1$, then 
 \[
\langle \mu + \rho, \epsilon_1-\epsilon_{j} \rangle \leq  \langle \mu + \rho, \epsilon_1-\epsilon_{j_1-1} \rangle <p \leq n_{1j}p,
 \]
 and if $j \geq j_1$, then
\begin{align*}
\langle \mu + \rho, \epsilon_1-\epsilon_{j} \rangle &= 
	\langle \mu + \rho, \epsilon_1-\epsilon_{j_1} \rangle + \langle \mu + \rho, \epsilon_{j_1}-\epsilon_{j} \rangle \\
				&= p + \langle \nu + \rho, \epsilon_{j_1}-\epsilon_{j} \rangle\\
				 &< (n_{j_1j}+1)p \\
				 & \leq n_{1j}p.
\end{align*}
The inequality $n_{j_1j} +1 \leq n_{1j}$ holds because 
\[
\langle \lambda + \rho, \epsilon_{1}-\epsilon_{j}\rangle \geq p + \langle \lambda + \rho, \epsilon_{j_1}-\epsilon_{j}\rangle
		\geq  p + (n_{j_1j}-1)p = n_{j_1j}p,
\]
where we set $n_{j_1j}=1$ if $j_1=j$. 
Therefore,  the third condition is also satisfied by $\mu$ when $i_1=1$, so now assume
$i_1 \geq 2$.

The problem reduces to finding positive integers $a_1, \dots, a_{i_1-1}$ such that 
(2) and (3) are completely satisfied by $\mu$.
Let us first assume that
\[ a = a_1 = \cdots = a_{i_1-1},\]
where $a > 0$. 
If $a < b_{j_1-2}/(i_1-1)$, then 
\begin{equation}\label{eqn:condition-1}
\begin{aligned}
\langle \mu + \rho, \epsilon_{1}-\epsilon_{j_1-1}\rangle &= \langle \mu + \rho, \epsilon_1-\epsilon_{i_1} \rangle + 
			\langle \mu + \rho, \epsilon_{i_1} - \epsilon_{j_1-1}\rangle \\ 
					&= (i_1-1)a + p-b_{j_1-2}  \\
					&< p,
\end{aligned}
\end{equation}
and hence $\mu$ satisfies $(2)$. 

It remains to determine sufficient conditions on $a$ so that $\mu$ satisfies (3). 
If $i > i_1$, then 
\[
\langle \mu + \rho, \epsilon_{i}-\epsilon_{j}\rangle = \langle \nu + \rho, \epsilon_{i}-\epsilon_{j}\rangle 
< n_{ij}p,
\]
by condition (3) on $\nu$.
Furthermore,  if $j < j_1$, then by \eqref{eqn:condition-1}, 
\begin{align*}
\langle \mu + \rho, \epsilon_{i}-\epsilon_{j}\rangle &\leq \langle \mu + \rho, \epsilon_{1}-\epsilon_{j_1-1}\rangle \\
	&=(i_{1}-1)a + p-(b_{j+1}+\cdots + b_{j_1-2}) \\
	&< p \leq n_{ij}p.
\end{align*}
Thus, \eqref{eqn:reduced-problem} is met for any $\epsilon_i - \epsilon_j$ with $i > i_1$ or $j < j_1$.

If $\epsilon_i - \epsilon_j \in \Phi^+$ is such that  $i \leq i_1$ and $j \geq j_1$, then 
$\epsilon_i-\epsilon_j = (\epsilon_i-\epsilon_{j_1}) + (\epsilon_{j_1}-\epsilon_j)$. 
Now  since $\lambda \in \mathcal{C}$, 
\[
\langle \lambda + \rho, \epsilon_i - \epsilon_{j_1}\rangle  \geq 
\langle \lambda + \rho, \epsilon_{i_1}-\epsilon_{j_1}\rangle > p,
\]
 and hence,
\begin{align*}
\langle \lambda + \rho, \epsilon_{i}-\epsilon_{j}\rangle &\geq p + \langle \lambda + \rho, \epsilon_{j_1}-\epsilon_{j}\rangle
		\geq  p + (n_{j_1j}-1)p = n_{j_1j}p.
\end{align*}
Thus, $n_{ij} \geq n_{j_1j}+1$, where we set $n_{j_1j}=1$ if $j_1=j$.

Also,
\begin{align*}
\langle \mu + \rho, \epsilon_{i}-\epsilon_{j}\rangle &= \langle \mu + \rho, \epsilon_{i}-\epsilon_{i_1}\rangle + 
  \langle \mu + \rho, \epsilon_{i_1}-\epsilon_{j_1}\rangle + \langle \mu + \rho, \epsilon_{j_1}-\epsilon_{j}\rangle \\
		&= (i_1-i)a+ p + \langle \nu + \rho, \epsilon_{j_1}-\epsilon_{j}\rangle,
\end{align*}
since $\langle \mu + \rho, \epsilon_{j_1}-\epsilon_{j}\rangle = \langle \nu + \rho, \epsilon_{j_1}-\epsilon_{j}\rangle$.
Hence, if $a$ is chosen so that for each $\epsilon_i - \epsilon_j \in \Phi^+$ with $i\leq i_1$ and $j\geq j_1$,  the
inequality
\begin{equation}\label{eqn:condition-2}
(i_1-i)a +\langle \nu + \rho, \epsilon_{j_1}-\epsilon_{j}\rangle  < n_{j_1j}p
\end{equation}
holds, 
then \eqref{eqn:reduced-problem} will hold for any $\epsilon_i - \epsilon_j$ with $i \leq i_1$ and $j \geq  j_1$ because 
\[
\langle \mu + \rho, \epsilon_{i}-\epsilon_{j}\rangle =  (i_1-i)a + \langle \nu + \rho, \epsilon_{j_1}-\epsilon_{j}\rangle + p < (n_{j_1j} + 1)p \leq n_{ij}p.
\]
Where a solution exists to \eqref{eqn:condition-2} for some $a > 0$, since
$\nu$ satisfies (3) by 
the inductive hypothesis, and hence 
$
\langle \nu + \rho, \epsilon_{j_1}-\epsilon_{j}\rangle < n_{j_1j}p.
$
Therefore, $\langle \mu + \rho, \epsilon_i - \epsilon_j \rangle < n_{ij}p$ for all $\epsilon_i-\epsilon_j \in \Phi^+$ which implies that 
$\mu$ satisfies (3). 

In summary, we shown that if 
\[
\mu + \rho = (a_1,\dots, a_n),
\]
where $a_k = b_k$ for $k\geq i_1$, $a_{i_1} = p-(b_{i_1+1}+\cdots +b_{j_1-1})$ and 
$a_k = a$ for $k=1,\dots, i_1-1$ such that $a > 0$ and satisfies \eqref{eqn:condition-1} and \eqref{eqn:condition-2}, then 
$\mu$ satisfies conditions (1), (2) and (3).  Therefore, the desired result follows by induction.

\end{proof}
Unfortunately, even when $p\geq n+1$ and $\lambda \in X(T)_+$,  
Lemma~\ref{lem:gamma-subroot} doesn't guarantee that
for each good positive subroot system $\Psi \subseteq \Gamma_{\lambda}$,
 there exists a weight 
$\mu \in X(T)_+$ such that $\Phi_{\mu,p}^+ \supseteq \Psi$ and 
 $C(\mu) \preceq C(\lambda)$. It only ensures the existence of a Euclidean point $\mu \in \mathcal{C}$ with the
 desired properties. 
   This issue will be clarified by the following lemma. 
  
    \begin{lem}\label{lem:latticept}
  Let $p\geq n+1$, then for every non-empty facette $F \subset \mathbb{E}$,
   $F \cap X(T) \neq\emptyset$.
   Equivalently, every non-empty facette contains a lattice point. 
  \end{lem}
  \begin{proof}
 The statement follows immediately from Proposition~\ref{prop:latticept}, since every subroot system of $\Phi$ is of the form 
 $w(\Phi_I)$ for some $I\subseteq \Delta$ and $w \in W$. 
   \end{proof}
  \begin{rem}
  The conclusion of the preceding lemma generally doesn't hold for other types. For example, let 
  $\Phi\subset \mathbb{R}^2$ be of type $C_2$ and let 
  $\alpha_1 = \epsilon_1 -\epsilon_2$ and $\alpha_2 = 2\epsilon_2$ denote the standard basis. 
  The corresponding fundamental weights are 
  $\omega_1 =  \epsilon_1$ and $\omega_2 = \epsilon_1 + \epsilon_2$, so that $\rho = 2\epsilon_1 + \epsilon_2$. 
  Consider the facette $F$, consisting of all $\lambda \in \overline{C_0}$ satisfying
  \begin{align*}
      \langle \lambda + \rho, \alpha_1^{\vee}\rangle &= 0 \\
  \langle \lambda + \rho, \alpha_1^{\vee} +  2\alpha_2^{\vee}\rangle &= p.
  \end{align*}
  If we write $\lambda + \rho = a_1\,\omega_1 + a_2\,\omega_2$, then 
  the first equation forces $a_1 = 0$, and so the second equation reduces to
   $2a_2 = p$. 
    Thus, the only integral solution occurs when $p$ is even and 
   $\lambda + \rho = (p/2)\,\omega_2$. 
  \end{rem}

\begin{prop}\label{prop:upperbound}
 Let $p \geq n+1$, then  for any $\lambda \in X(T)_+$, 
 $
 V_{G_1}(T(\lambda)) \subseteq \overline{\mathcal{O}_{s(\lambda)^t}}.
 $
 \end{prop}
 \begin{proof}
  By Lemma~\ref{lem:gamma-subroot}, 
    for each good positive subroot system $\Psi \subseteq \Gamma_{\lambda}$, there exists
     $\mu \in \mathcal{C}$ such that   $\Psi\subseteq \Phi_{\mu,p}^+$ and  $C(\mu) \preceq C(\lambda)$. 
      Thus, 
  $\pi(\Psi) \leq \pi(\Phi_{\mu,p}^+) = d(\mu)$.
     Furthermore, by Lemma~\ref{lem:latticept}, $F(\mu) \cap X(T)_+ \neq \emptyset$, so 
   we may assume that $\mu \in X(T)_+$. 
  Then
  \begin{align*}
  V_{G_1}(T(\lambda)) &\subseteq V_{G_1}(H^0(\mu)) & \text{(by Proposition~\ref{cor:ind-tilting})} \\
                                     &= \overline{\mathcal{O}_{d(\mu)^t}} &  \text{(by \cite[Theorem 6.2.1]{npv2002})} \\
                                     &\subseteq \overline{\mathcal{O}_{\pi(\Psi)^t}}.
  \end{align*}
  Hence,  for each good positive subroot system $\Psi\subseteq \Gamma_{\lambda}$,
  $V_{G_1}(T(\lambda))\subseteq \overline{\mathcal{O}_{\pi(\Psi)^t}}$.
  Therefore, by Proposition~\ref{cor:good-part},  
  $V_{G_1}(T(\lambda)) \subseteq \overline{\mathcal{O}_{s(\lambda)^t}}$
  since $s(\lambda)^t$ is the greatest partition satisfying $s(\lambda)^t \leq \pi(\Psi)^t$ for 
  all good $\Psi\subseteq \Gamma_{\lambda}$. 
%
%
\end{proof}

\section{Quantum groups}\label{sec:quantum}
\subsection{}
This section will follow the notation and conventions in \cite{bnpp} and \cite[Appendix H]{jan2003}.
Let $\mathfrak{g}=\mathfrak{g}_{\mathbb{C}}$ denote a finite-dimensional, complex, semisimple Lie algebra 
and let $G_{\mathbb{Z}}$ denote the split, semisimple, simply connected algebraic group scheme such that
$\mathfrak{g} = \text{Lie}(G)\otimes_{\mathbb{Z}}\mathbb{C}$. 
Denote by $\mathbb{U}_q(\mathfrak{g})$, the  \emph{quantum enveloping algebra} with indeterminate 
$q \in \mathbb{Q}(q)$ and generators $E_{\alpha}$, $F_{\alpha}$, $K_{\alpha}$ and 
$K^{-1}_{\alpha}$ for $\alpha \in \Delta$, satisfying the quantized Serre relations (\cite[H.2]{jan2003}). 

For $\mathscr{A} = \mathbb{Z}[q,q^{-1}]$, let 
$\mathbb{U}_q^{\mathscr{A}}(\mathfrak{g})$ 
be the \emph{Lusztig $\mathscr{A}$-form} of $\mathbb{U}_q(\mathfrak{g})$, which is the $\mathscr{A}$-subalgebra
generated by the divided powers $E_{\alpha}^{(m)}$, $F_{\alpha}^{(m)}$ and $K_{\alpha}^{\pm 1}$ 
(cf. \cite[H.5]{jan2003}).
If $\Gamma$ is an $\mathscr{A}$-algebra, set 
$\mathbb{U}_{\Gamma}(\mathfrak{g}) = \mathbb{U}_q^{\mathscr{A}}(\mathfrak{g})\otimes_{\mathbb{Z}}\mathbb{C}$. Finally, for any  
$\ell >1$ and any primitive $\ell^{th}$ root of unity $\zeta \in \mathbb{C}$,  give $\Gamma = \mathbb{C}$
the structure of an $\mathscr{A}$-algebra by sending $q \mapsto \zeta$, and write
$\mathbb{U}_{\zeta}(\mathfrak{g}) = \mathbb{U}_{\Gamma}(\mathfrak{g})$.

It is well known that the category of type 1 integrable representations for $\mathbb{U}_{\zeta}(\mathfrak{g})$ shares many properties
in common with the category of modular representations for $G_k$, where  
$k$ an algebraically closed field of characteristic $p>0$ (cf. \cite[Appendix H]{jan2003} for an overview).  
Thus, every $\mathbb{U}_{\zeta}(\mathfrak{g})$-module will be assumed to be type 1 and integrable.
Suppose $I \subseteq \Delta$ is a subset and $\mathfrak{l}_I$ and $\mathfrak{p}_I$ 
are the corresponding Levi and (negative) parabolic subalgebras of $\mathfrak{g}$, then  one 
can define \emph{Levi and parabolic subalgebras} 
 $\mathbb{U}_q(\mathfrak{l}_I)$ and
$\mathbb{U}_q(\mathfrak{p}_I)$ 
 of the quantum enveloping algebra $\mathbb{U}_q(\mathfrak{g})$ and, by specialization, the subalgebras
 $\mathbb{U}_{\zeta}(\mathfrak{l}_I)$ and
$\mathbb{U}_{\zeta}(\mathfrak{p}_I)$ of $\mathbb{U}_{\zeta}(\mathfrak{g})$ (cf. \cite[2.5]{bnpp}).
It is possible to define the \emph{induction functor},
 \[
 \myind{\mathbb{U}_{\zeta}(\mathfrak{p})}{\mathbb{U}_{\zeta}(\mathfrak{g})}(M) = H^0(\mathbb{U}_{\zeta}(\mathfrak{g})/\mathbb{U}_{\zeta}(\mathfrak{p}), M)
 \] 
for any $\mathbb{U}_{\zeta}(\mathfrak{p})$-module $M$ (\cite[2.4]{apw}). When dealing with the \emph{Borel subalgebra} $\mathbb{U}_{\zeta}(\mathfrak{b})$, we will write
$H^0_{\zeta}(M) = \myind{\mathbb{U}_{\zeta}(\mathfrak{b})}{\mathbb{U}_{\zeta}(\mathfrak{g})}(M)$ for any $\mathbb{U}_{\zeta}(\mathfrak{b})$-module $M$. 

Let $X$ denote the weight lattice for $\mathbb{U}_{\zeta}(\mathfrak{g})$ and let $X^+$ denote the cone of dominant weights. 
 For each $\lambda \in X^+$,  $L_{\zeta}(\lambda)=\text{soc}\,(H^0_{\zeta}(\lambda))$ 
 is the corresponding simple highest weight
 module (cf. \cite[Corollary 6.2]{apw}).
 The Weyl modules are defined by $V_{\zeta}(\lambda)=H_{\zeta}^0(-w_0\lambda)^*$.
The tilting modules are defined  in the same way as for algebraic groups. It follows that 
 for each weight $\lambda \in X^+$, there exists a unique
indecomposable tilting module $T_{\zeta}(\lambda)$
for $\mathbb{U}_{\zeta}(\mathfrak{g})$ (\cite[H.15]{jan2003}). It was proven by Soergel that under some slight restrictions on $\ell$, the formal characters of these modules
are determined by certain parabolic Kazhdan-Lusztig polynomials (cf. \cite{soergel1} and \cite{soergel2}).

Much work has been done in studying the cohomology of the finite-dimensional
Hopf algebra $\mathfrak{u}_{\zeta}(\mathfrak{g}) \unlhd \mathbb{U}_{\zeta}$, known as the small quantum group 
(\cite[2.2]{bnpp}). 
For instance, in \cite[Theorem 3]{gk} it was shown that when $\ell > h$, 
\[
\text{H}^{\text{ev}}(\mathfrak{u}_{\zeta}(\mathfrak{g}), \mathbb{C})_{\text{red}} \cong \mathbb{C}[\mathcal{N}].
\]
Thus, $\text{max-Spec}(\text{H}^{\text{ev}}(\mathfrak{u}_{\zeta}(\mathfrak{g}), \mathbb{C})) = \mathcal{N}$,
where $\mathcal{N}\subseteq \mathfrak{g}$ is the nilpotent cone of $\mathfrak{g}$.  
To each $u_{\zeta}(\mathfrak{g})$-module $M$, there exists a support variety 
$V_{\mathfrak{u}_{\zeta}(\mathfrak{g})}(M)\subseteq \mathcal{N}$. If $M$ has the structure of a 
$\mathbb{U}_{\zeta}(\mathfrak{g})$-module, then 
$V_{\mathfrak{u}_{\zeta}(\mathfrak{g})}(M)$ is in fact a $G_{\mathbb{C}}$-stable subvariety of $\mathcal{N}$
(cf. \cite[8.1]{bnpp}).

\subsection{}
Recalling Definition~\ref{def:tilting-category}, let
$\mathcal{T}_{\zeta}$ denote the full subcategory of all finite-dimensional 
 tilting modules for $\mathbb{U}_{\zeta}(\mathfrak{g})$.
The thick tensor ideals of $\mathcal{T}_{\zeta}$ have been classified by Ostrik (cf. \cite[Theorem 4.5]{ostrik97}). 
More specifically, it was shown that 
$\langle T_{\zeta}(\mu) \rangle = \langle T_{\zeta}(\lambda) \rangle$ 
if and only if $\lambda, \mu \in X^+$ lie in the same weight cell.
In further analogy with the algebraic group case, there is a connection between thick tensor ideals
and support varieties for quantum tilting modules.
 \begin{lem}\label{lem:sup-ideal}
 Let $M, N \in \mathcal{T}_{\zeta}$ be tilting modules with
  $\langle M \rangle \subseteq \langle N \rangle$, then 
 $V_{\mathfrak{u}_{\zeta}(\mathfrak{g})}(M) \subseteq V_{\mathfrak{u}_{\zeta}(\mathfrak{g})}(N)$.
 \end{lem}
 \begin{proof}
 Since $\langle M \rangle \subseteq \langle N \rangle$, then by definition 
 there exists some $L \in \mathcal{T}_{\zeta}$ such that $M \mid N \otimes L$,
  and hence 
   $V_{\mathfrak{u}_{\zeta}(\mathfrak{g})}(M) \subseteq V_{\mathfrak{u}_{\zeta}(\mathfrak{g})}(N\otimes L) \subseteq 
   V_{\mathfrak{u}_{\zeta}(\mathfrak{g})}(N).$
 \end{proof}

In Section~\ref{sec:weight-cells}, it was stated that the varieties 
$V_{\mathfrak{u}_{\zeta}(\mathfrak{g})}(T_{\zeta}(\lambda))$ have been computed 
 for all types when $\ell > h$ by Ostrik and Bezrukavnikov.
 For convenience, we shall state here what was proven.
\begin{thm}\label{thm:hum-quantum}
 Let $\mathfrak{g}$ be a complex semisimple Lie algebra, and let
 $\zeta \in \mathbb{C}$ be a primitive $\ell^{th}$ root of unity with $\ell>h$, odd (and not divisible by 3 if $\mathfrak{g}$ has a component of 
 type $G_2$).  For each $w \in W^+_{\ell}$, let $c_{[w]} \subseteq \mathcal{C}$ be the corresponding weight cell, and let 
 $\mathcal{O}_{[w]}$ denote the orbit associated to $c_{[w]}$ by the Lusztig bijection. Then
 if $\lambda \in c_{[w]}\cap X^+,$
\[
V_{\mathfrak{u}_{\zeta}(\mathfrak{g})}(T_{\zeta}(\lambda)) = \overline{\mathcal{O}_{[w]}}.
\]
\end{thm}

 For any tilting module $M$ for $G$, it is well known that $M|_{[L_I.L_I]}$ is a tilting module whenever $I\subseteq \Delta$,
 $L_I$ is a Levi-factor for $G$ and $[L_I,L_I]$ is the derived subgroup of $L_I$ (cf. \cite[Proposition II.4.24]{jan2003}).
An analogous result also holds for quantum groups.
 \begin{prop}\label{prop:nt-prop}
Let $\mathfrak{g}$ be a complex semisimple
 Lie algebra, and let
 $\zeta \in \mathbb{C}$ be a primitive $\ell^{th}$ root of unity where $\ell$ is odd (and not divisible by
$3$ if $\mathfrak{g}$ has a component of type $G_2$) and is such that $\langle \omega_i + \rho, \alpha_0^{\vee}\rangle < \ell$ for all
fundamental weights $\omega_1,\dots, \omega_n$. 
Then for each $I\subseteq \Delta$,  $\mathbb{U}_{\zeta}([\mathfrak{l}_I, \mathfrak{l}_I])$ is the Hopf subalgebra of $\mathbb{U}_{\zeta}(\mathfrak{g})$ generated by $E_{\alpha}^{(m)}$, $F_{\alpha}^{(m)}$, 
$K_{\alpha}^{\pm 1}$ for $\alpha \in I$, and 
for any $\mathbb{U}_{\zeta}(\mathfrak{g})$ tilting module $M$, the restricted module $M|_{\mathbb{U}_{\zeta}([\mathfrak{l}_I,\mathfrak{l}_I])}$ is a $\mathbb{U}_{\zeta}([\mathfrak{l}_I,\mathfrak{l}_I])$ tilting 
module.
\end{prop}
\begin{proof}
Begin by observing that for each fundamental weight $\omega \in X^+$, the restriction
 $H^0_{\zeta}(\omega)|_{\mathbb{U}_{\zeta}([\mathfrak{l}_I,\mathfrak{l}_I])}$ is a tilting module. This is due to the fact that all of the weights of 
$H^0_{\zeta}(\omega)|_{\mathbb{U}_{\zeta}([\mathfrak{l}_I,\mathfrak{l}_I])}$ are $\ell$-miniscule (i.e. they satisfy $\langle \nu + \rho, \alpha_0^{\vee}\rangle < \ell$),
 and hence the restricted module must be a semisimple tilting module. 
The proposition follows by adapting the argument in  \cite[Proposition 3.1]{nt} to the quantum setting. 
\end{proof}
\begin{rem}
If $\mathfrak{g} = \mathfrak{sl}_{n+1}(\mathbb{C})$, then the condition that $\langle \omega_i + \rho, \alpha_0^{\vee}\rangle < \ell$ for all
fundamental weights $\omega_1,\dots, \omega_n$, is satisfied precisely when $\ell > n+1$. 
\end{rem}

\subsection{}
Suppose now that $\zeta \in \mathbb{C}$ is a $p^{th}$ root of unity, where $p$ is a prime number, and let 
$k$ be an algebraically closed field of characteristic $p$. Let
 $\mathcal{T}$ denote the subcategory of tilting modules
for the algebraic group $G= G_k$, and identify $X=X(T)$, where $X$ is the weight lattice for 
$\mathbb{U}_{\zeta}(\mathfrak{g})$.  By \cite[5.3]{and1998}, for each tilting module $M$ of $G$, there exists 
a quantum tilting module, denoted by $M_{\zeta}$, for $\mathbb{U}_{\zeta}(\mathfrak{g})$ 
satisfying $\text{ch}\,(M_{\zeta}) = \text{ch}\,(M)$.
More specifically, if $\lambda \in X(T)_+$ is arbitrary and $M= T(\lambda)$, then
\begin{equation*}
T(\lambda)_{\zeta} = T_{\zeta}(\lambda) \oplus  \bigoplus_{\mu \uparrow \lambda, \mu \neq \lambda}a_{\mu}T_{\zeta}(\mu).
\end{equation*}
In particular, since $\text{ch}\, (T(\lambda)_{\zeta}) = \text{ch}\, (T(\lambda))$, then 
\begin{equation}\label{eqn:lifting}
\text{ch}\,(T(\lambda)) = \text{ch}\,(T_{\zeta}(\lambda)) +  \sum_{\mu \uparrow \lambda, \mu \neq \lambda}a_{\mu}\text{ch}\, (T_{\zeta}(\mu)).
\end{equation}

 From the standard properties of the assignment $V_{\mathfrak{u}_{\zeta}(\mathfrak{g})}(-)$, one gets
 $V_{\mathfrak{u}_{\zeta}(\mathfrak{g})}(T_{\zeta}(\lambda)) \subseteq V_{\mathfrak{u}_{\zeta}(\mathfrak{g})}(T(\lambda)_{\zeta})$ (cf. \cite[Lemma 3.4]{ostrik98}).
 This gives us an immediate corollary to Theorem~\ref{thm:hum-quantum}. 
 \begin{cor}\label{cor:lift}
  Let $G$ be a semisimple, simply connected  
 algebraic group over a field $k$ of characteristic $p >h$, and let 
 $\mathbb{U}_{\zeta}(\mathfrak{g})$ be the corresponding quantum group, where $\zeta \in \mathbb{C}$ is a primitive $p^{th}$ root 
 of unity. 
 Then for each $w \in W_p^+$, and
  $\lambda \in  c_{[w]}\cap X(T)_+,$
\[
V_{\mathfrak{u}_{\zeta}(\mathfrak{g})}(T(\lambda)_{\zeta}) \supseteq \overline{\mathcal{O}_{[w]}}.
\]
 \end{cor}
 \subsection{}
An interesting problem would be to understand how the support varieties for 
 tilting modules of the form $T(\lambda)$ and $T(\lambda)_{\zeta}$, with $\lambda \in X(T)_+$, are related. 
 It is well known that when the characteristic $p$ is \emph{good} (in particular if $p >h$), the classification 
 and structure of the $G_{\mathbb{C}}$ orbits on the complex nilpotent cone 
 $\mathcal{N}_{\mathbb{C}}=\mathcal{N}(G_{\mathbb{C}})$ coincide
 with the $G_{k}$ orbits on $\mathcal{N}_k=\mathcal{N}(G_k)$ (cf. \cite{cm} for the complex case and \cite{pom1977} for 
 the positive characteristic case). This implies that each orbit $\mathcal{O}_{\mathbb{C}}$ in $\mathcal{N}_{\mathbb{C}}$ uniquely
 corresponds to an orbit $\mathcal{O}_k$ in $\mathcal{N}_k$. Moreover, if 
 \[
 \overline{\mathcal{O}_{\mathbb{C}}} = \mathcal{O}^1_{\mathbb{C}}\cup \cdots \cup \mathcal{O}^m_{\mathbb{C}}
 \]
 for some orbits $\mathcal{O}^1_{\mathbb{C}},\dots, \mathcal{O}^m_{\mathbb{C}}$, then 
 \[
 \overline{\mathcal{O}_{k}} = \mathcal{O}^1_{k}\cup \cdots \cup \mathcal{O}^m_{k}.
 \]
It follows that any $G_{\mathbb{C}}$-stable closed subvariety $V_{\mathbb{C}}\subseteq \mathcal{N}_{\mathbb{C}}$ uniquely corresponds to 
a $G_k$-stable closed subvariety $V_{k}\subseteq \mathcal{N}_{k}$. 
We now state an interesting conjecture which would realize this correspondence by taking support varieties of tilting modules.
\begin{conj}\label{conj:correspondence}
 Let $G$ be a semisimple, simply connected  
 algebraic group over a field $k$ of characteristic $p >h$, and let 
 $\mathbb{U}_{\zeta}(\mathfrak{g})$ be the corresponding quantum group, where $\zeta \in \mathbb{C}$ is a primitive $p^{th}$ root 
 of unity. 
  Then for any tilting module $M$ for $G_k$,  
  $V_{\mathfrak{u}_{\zeta}(\mathfrak{g})}(M_{\zeta}) = V_{\mathbb{C}}$ if and only if
 $V_{G_1}(M) = V_{k}$, where $V_k$ is the unique $G_k$-stable subvariety of $\mathcal{N}_k$ corresponding to 
 $V_{\mathbb{C}}$. 
 \end{conj}
 \begin{rem}\label{rem:correspondence}
 The truth of this conjecture would imply that the correspondence between the $G_{\mathbb{C}}$
 and $G_k$-stable closed subvarieties of $\mathcal{N}_{\mathbb{C}}$ and $\mathcal{N}_k$ described above, can be 
 established by taking support varieties of tilting modules. 
 In fact, if $p>h$, then the conjecture will follow if both Conjecture~\ref{conj:hum2} holds and an
 analogous conjecture holds for tilting modules of the form $T(\lambda)_{\zeta}$, where $\lambda \in X(T)_+$.
 \end{rem}

The following lemma verifies this conjecture for the trivial orbit closures $\{0\}_{\mathbb{C}}\subseteq \mathcal{N}_{\mathbb{C}}$
 and $\{0\}_k\subseteq \mathcal{N}_k$.
 \begin{lem}\label{lem:proj-lift}
 Let $G$ be a semisimple, simply connected  
 algebraic group over a field $k$ of characteristic $p >h$ and let 
 $\mathbb{U}_{\zeta}(\mathfrak{g})$ be the corresponding quantum group, where $\zeta \in \mathbb{C}$ is a primitive $p^{th}$ root 
 of unity.  Then a tilting $G$-module $M$ is $G_1$-projective if and only if
 $M_{\zeta}$ is $\mathfrak{u}_{\zeta}(\mathfrak{g})$-projective. 
\end{lem}
\begin{proof}
Without loss of generality, we may assume that $M = T(\lambda)$ for some $\lambda \in X(T)_+$. By 
\cite[Lemma E.8]{jan2003}, it follows that $T(\lambda)$ is $G_1$-projective if and only if 
$\langle \lambda, \alpha^{\vee}\rangle \geq p-1$
for all $\alpha \in \Delta$.  The analogous statement also holds in the quantum setting for $T_{\zeta}(\lambda)$. 
Since $$T(\lambda)_{\zeta} = T(\lambda) \oplus \bigoplus_{\mu \uparrow \lambda, \mu \neq \lambda} a_{\mu}T_{\zeta}(\mu),$$ then
$T(\lambda)_{\zeta}$ is projective if and only if $a_{\mu} = 0$ for any 
$\mu \in X(T)_+$ satisfying $\langle \mu, \alpha^{\vee}\rangle < p-1$ for some $\alpha \in \Delta$.

Now observe that  $T(\lambda)$ is $G_1$-projective if and only if 
$\lambda = (p-1)\rho + \nu$,
where $\nu \in X(T)_+$. 
Since $T((p-1)\rho) = L((p-1)\rho)$ is a simple tilting module, then
$
T_{\zeta}((p-1)\rho) = T((p-1)\rho)_{\zeta}
$ 
because both modules are equal to $L_{\zeta}((p-1)\rho)$. 
By highest weight considerations, 
$$T(\lambda)_{\zeta} \mid T_{\zeta}((p-1)\rho) \otimes T(\nu)_{\zeta},$$
and so $T(\lambda)_{\zeta}$ is projective if  $T(\lambda)$ is $G_1$-projective. Likewise, if $T(\lambda)_{\zeta}$ is 
projective, then  $T_{\zeta}(\lambda)$ is projective, and thus $\lambda = (p-1)\rho + \nu$
for some $\nu \in X(T)_+$, which implies that $T(\lambda)$ is $G_1$-projective. 
\end{proof}

Conjecture~\ref{conj:correspondence} can also be verified for the principal orbit closures:
$\mathcal{N}_{\mathbb{C}}$ and $\mathcal{N}_k$. 
\begin{lem}\label{lem:principal-lift}
 Let $G$ be a semisimple, simply connected  
 algebraic group over a field $k$ of characteristic $p >h$ and let 
 $\mathbb{U}_{\zeta}(\mathfrak{g})$ be the corresponding quantum group, where $\zeta \in \mathbb{C}$ is a primitive $p^{th}$ root 
 of unity.  Then any tilting $G$-module $M$ satisfies $V_{G_1}(M) = \mathcal{N}_k$ if and only if 
 $V_{\mathfrak{u}_{\zeta}(\mathfrak{g})}(M_{\zeta}) = \mathcal{N}_{\mathbb{C}}$. 
\end{lem}
\begin{proof}
Let $M = \sum_{\lambda \in X(T)_+} a_{\lambda}\, T(\lambda)$ be an arbitrary tilting module for $G$, then 
\[
M_{\zeta} = \sum_{\lambda \in X(T)_+} a_{\lambda}\, T(\lambda)_{\zeta} = \sum_{\lambda \in X(T)_+} b_{\lambda}\, T_{\zeta}(\lambda).
\]
By using translation identities, it can be deduced that $V_{G_1}(M) = \mathcal{N}_k$ if and only if $a_{\lambda} > 0$ for some $\lambda \in C_0$
(cf. \cite[Proposition E.11]{jan2003}).
By the same argument, 
$V_{\mathfrak{u}_{\zeta}(\mathfrak{g})}(M_{\zeta}) = \mathcal{N}_{\mathbb{C}}$ if and only if $b_{\lambda} > 0$ for some $\lambda \in C_0$. 
On the other hand, by \cite[Proposition E.12]{jan2003}, it follows that for each $\lambda \in C_0$,
\[
a_{\lambda} = \sum_{w \in W_p^+}(-1)^{\ell(w)}[M:\mych H^0(w\cdot \lambda)].
\]
Moreover, since $\mych H^0_{\zeta}(\mu) = \mych H^0(\mu)$ for any $\mu \in X(T)_+$ and since $\mych M = \mych M_{\zeta}$,
 then by the same argument,
\[
b_{\lambda} = \sum_{w \in W_p^+}(-1)^{\ell(w)}[M:\mych H^0(w\cdot \lambda)]
\]
for any $\lambda \in C_0$. 
Thus, $a_{\lambda} = b_{\lambda}$ for each $\lambda \in C_0$, and hence  $V_{G_1}(M) = \mathcal{N}_k$ if and only if 
$V_{\mathfrak{u}_{\zeta}(\mathfrak{g})}(M_{\zeta}) = \mathcal{N}_{\mathbb{C}}$. 

\end{proof}

\subsection{An interesting result in type $A_n$}
Let $G=SL_{n+1}(k)$ with $p\geq n+1$,  let $\mathfrak{g} = \mathfrak{sl}_{n+1}(\mathbb{C})$, and let 
$\zeta \in \mathbb{C}$ be a primitive $p^{th}$ root of unity. By 
Corollary~\ref{cor:lift},
$V_{\mathfrak{u}_{\zeta}(\mathfrak{g})}(T(\lambda)_{\zeta}) \supseteq \overline{\mathcal{O}_{s(\lambda)^t}}$ for any 
$\lambda \in X(T)_+$. 
It can also be verified  that Proposition~\ref{cor:ind-tilting} holds in the quantum setting. 
Therefore, the proof of Proposition~\ref{prop:upperbound} may
be adapted to the quantum group setting, to yield the following proposition. 
\begin{prop}\label{prop:quantum-upper}
Let $p > n+1$, then  for each $\lambda \in X(T)_+$, 
 $
 V_{\mathfrak{u}_{\zeta}}(T(\lambda)_{\zeta}) =\overline{\mathcal{O}_{s(\lambda)^t}}.
 $
\end{prop}

By Remark~\ref{rem:correspondence}, this proposition can be combined with Theorem~\ref{thm:main} to 
prove
Conjecture~\ref{conj:correspondence} in the type $A$ case.

\section{The lower bound}\label{sec:type-A}
\subsection{}\label{subsection:lower-setup}
Let $G = SL_{n+1}(k)$,  where $k$ is an algebraically closed field of characteristic $p > 0$, 
let $\mathfrak{g} = \mathfrak{sl}_{n+1}(\mathbb{C})$ and let 
$\zeta \in \mathbb{C}$ be a primitive $p^{th}$ root of unity.
 The goal of this section will be to show that for each $\lambda \in X(T)_+$,
\begin{equation}\label{eq:hum}
V_{G_1}(T(\lambda)) \supseteq \overline{\mathcal{O}_{s(\lambda)^t}}.
\end{equation}
For any partition $\pi= (p_1,p_2,\dots, p_r )\in \mathcal{P}$, we define the subgroup scheme
\[
SL_{\pi} = SL_{p_1} \times SL_{p_2} \times \cdots \times SL_{p_r} \subseteq SL_{n+1}, 
\]
and let  $\mathfrak{sl}_{\pi} = \text{Lie}(SL_{\pi})$ denote its Lie algebra.  
If $\alpha_i = \epsilon_i-\epsilon_{i+1}$ for $i=1,\dots, n$ and
\begin{equation*}
I_{\pi} = \{\alpha_1,\dots, \alpha_{p_1-1}\}\cup \{\alpha_{p_1+1},\dots, \alpha_{p_1+p_2-1}\}\cup \cdots 
          \cup \{\alpha_{p_1+\cdots  p_{r-1}+1}, \dots, \alpha_{p_1+\cdots + p_r-1}\} \subseteq \Delta,
\end{equation*}
then $SL_{\pi} = [L_{I_{\pi}},L_{I_{\pi}}]$ is the derived subgroup of the corresponding Levi factor $L_{I_{\pi}}$. 
 For notational simplicity, we will set $H_{\pi}= SL_{\pi}(k)$ and
 $\mathfrak{h}_{\pi}=\mathfrak{sl}_{\pi}(\mathbb{C})$.

In Proposition~\ref{prop:nt-prop}, it was shown that, as in the algebraic setting, there is a natural
inclusion of quantum groups $\mathbb{U}_{\zeta}(\mathfrak{h}_{\pi}) \hookrightarrow \mathbb{U}_{\zeta}(\mathfrak{g})$ such that tilting modules 
for $\mathbb{U}_{\zeta}(\mathfrak{g})$ restrict to tilting modules for $\mathbb{U}_{\zeta}(\mathfrak{h}_{\pi})$. 
   The following lemma is a well known fact about nilpotent orbits in type $A$ (cf. \cite[Theorem 8.2.14]{cm}). 
 \begin{lem}\label{lem:a_n-levi}
 For any partition $\pi \in \mathcal{P}$, let  
 $x_{\pi}$ denote the nilpotent matrix which is a direct sum of Jordan blocks whose sizes are given by the parts of $\pi$,
 then $x_{\pi} \in \mathfrak{sl}_{\pi}(k)$ and the orbit $H_{\pi}\cdot x_{\pi}$ is dense in $\mathcal{N}(H_{\pi})$.
Moreover, since $\mathcal{O}_{\pi} = G\cdot x_{\pi}$, then $\overline{G\cdot \mathcal{N}(H_{\pi})} = \overline{\mathcal{O}_{\pi}}$.
%
\end{lem}

  By the naturality of support varieties, we can identify
\begin{equation}\label{eq:sl2-natural}
V_{(H_{\pi})_1}(M|_{H_{\pi}}) = V_{G_1}(M) \cap \mathfrak{sl}_{\pi}(k) 
\end{equation}
for each $G$-module $M$. 
Under this identification, $x_{\pi} \in V_{G_1}(M)$ if and only if 
$V_{(H_{\pi})_1}(M|_{H_{\pi}}) = \mathcal{N}(H_{\pi})$.

 \subsection{}
  We now have enough to  proceed with a proof of the lower bound. But first, 
  for notational convenience,  the following terminology will be introduced. 
 \begin{defn}
 A module $M$ is said to have \emph{full support}, if its support variety is maximal. 
 For instance, if $M$ is a $G_1$-module with $p>h$, then $M$ has full support provided
 $V_{G_1}(M) = \mathcal{N}(G)$.  
 \end{defn}
 As mentioned in the introduction,  the following proposition, along with Proposition~\ref{prop:upperbound} may be combined to give
 Theorem~\ref{thm:main}.
 \begin{prop}\label{prop:lowerbound}
 Let $p > n+1$, then  for any $\lambda \in X(T)_+$, 
 $
 V_{G_1}(T(\lambda)) \supseteq \overline{\mathcal{O}_{s(\lambda)^t}}.
 $
 \end{prop}
 \begin{proof}
For any partition $\pi \in \mathcal{P}$ and $\lambda \in X(T)_+$ satisfying $s(\lambda)^t = \pi$, it follows from 
 \eqref{eq:sl2-natural} that  this proposition will hold if
$T(\lambda)|_{H_{\pi}}$ has full support.
By Proposition~\ref{prop:nt-prop}, the module  $(T(\lambda)_{\zeta})|_{\mathbb{U}_{\zeta}(\mathfrak{h}_{\pi})}$ is a quantum tilting module, 
and thus $(T(\lambda)|_{H_{\pi}})_{\zeta} = (T(\lambda)_{\zeta})|_{\mathbb{U}_{\zeta}(\mathfrak{h}_{\pi})}$. So by Lemma~\ref{lem:principal-lift}, $T(\lambda)|_{H_{\pi}}$ will have full support if and only if $(T(\lambda)_{\zeta})|_{\mathbb{U}_{\zeta}(\mathfrak{h}_{\pi})}$
has full support.
Since $$T_{\zeta}(\lambda) \mid T(\lambda)_{\zeta},$$ then
$V_{\mathfrak{u}_{\zeta}(\mathfrak{g})}(T_{\zeta}(\lambda)) \subseteq V_{\mathfrak{u}_{\zeta}(\mathfrak{g})}(T(\lambda)_{\zeta})$,
and hence, by \cite{ostrik97}, \cite[Lemma 6.4]{ostrik98} and \cite[Theorem 6.8]{ostrik98}, 
 there exist $\mu, \nu  \in X(T)_+$ such that $T(\mu) = H^0(\mu)$, $s(\mu) = s(\lambda)$ and 
$T_{\zeta}(\mu) \mid T(\lambda)_{\zeta}\otimes T_{\zeta}(\nu).$
Moreover,  $T_{\zeta}(\mu) = T(\mu)_{\zeta}$ and $T_{\zeta}(\nu) \mid T(\nu)_{\zeta}$ together give 
\begin{equation}\label{eqn:weyl-divis}
 T(\mu)_{\zeta} \oplus M = (T(\lambda) \otimes T(\nu))_{\zeta},
\end{equation}
for some $\mathbb{U}_{\zeta}(\mathfrak{g})$ tilting module $M$. Therefore, by Proposition~\ref{prop:nt-prop}
\begin{equation}\label{eqn:weyl-divis2}
 T(\mu)_{\zeta}|_{\mathbb{U}_{\zeta}(\mathfrak{h}_{\pi})} \oplus M|_{\mathbb{U}_{\zeta}(\mathfrak{h}_{\pi})} = N_{\zeta},
\end{equation}
where $N = (T(\lambda) \otimes T(\nu))|_{H_{\pi}}$. Since $T(\mu)=H^0(\mu)$, then by 
\cite[Theorem 6.2.1]{npv2002}, $V_{G_1}(T(\mu)) = \overline{\mathcal{O}_{\pi}}$ and so \eqref{eq:sl2-natural},
implies that
$V_{(H_{\pi})_1}(T(\mu)|_{H_{\pi}}) = \mathcal{N}(H_{\pi})$. Also,  by \eqref{eqn:weyl-divis2} and Lemma~\ref{lem:principal-lift}, 
 $(T(\mu)|_{H_{\pi}})_{\zeta} = (T(\mu)_{\zeta})|_{\mathbb{U}_{\zeta}(\mathfrak{h}_{\pi})}$, and hence  $N$ and $N_{\zeta}$ have full support. 
Thus, $x_{\pi} \in V_{G_1}(T(\lambda)\otimes T(\nu)) \subseteq V_{G_1}(T(\lambda))$, where the inclusion follows from the tensor product identity for support varieties.
\end{proof}

\end{document}